\documentclass[10pt]{amsart}
\usepackage{latexsym}
\usepackage{amsmath}
\usepackage{amssymb}
\usepackage{mathrsfs}
\usepackage{graphicx}
\usepackage{color}
\usepackage{pgfpages}
\usepackage{ifthen}
\usepackage{leftidx,tensor}
\usepackage[T1]{fontenc}
\usepackage[latin1]{inputenc}
\usepackage{mathtools}
\usepackage{comment}
\usepackage{dsfont}
\usepackage[nocompress]{cite}
\newcommand{\citep}{\cite}

\usepackage[shortlabels]{enumitem}
\usepackage{aliascnt}
\usepackage[bookmarks=true,pdfstartview=FitH, pdfborder={0 0 0}, colorlinks=true,allcolors=black,linkcolor=blue,citecolor=red]{hyperref}
\usepackage{bbm}
\usepackage{nicefrac}

\theoremstyle{plain}
\newtheorem{thm}{Theorem}[section]
\newaliascnt{cor}{thm}
\newaliascnt{prop}{thm}
\newaliascnt{lem}{thm}
\newtheorem{cor}[cor]{Corollary}
\newtheorem{prop}[prop]{Proposition}
\newtheorem{lem}[lem]{Lemma}
\aliascntresetthe{cor}
\aliascntresetthe{prop}
\aliascntresetthe{lem} 
\theoremstyle{definition}
\newaliascnt{defn}{thm}
\newaliascnt{asu}{thm}
\newaliascnt{con}{thm}

\newtheorem{asu}[asu]{Assumption}

\aliascntresetthe{defn}
\aliascntresetthe{asu}
\aliascntresetthe{con}
\newcounter{stp}
\newcounter{stpi}
\newcounter{stpci}
\newcounter{stpiii}

\newtheorem{step}[stp]{Step}

\theoremstyle{plain}
\newaliascnt{rem}{thm}
\newaliascnt{exa}{thm}
\newaliascnt{masu}{thm}
\newaliascnt{nota}{thm}
\newaliascnt{sett}{thm}
\newtheorem{rem}[rem]{Remark}

\aliascntresetthe{rem}
\aliascntresetthe{exa}
\aliascntresetthe{masu}
\aliascntresetthe{nota}
\aliascntresetthe{sett}
\numberwithin{equation}{section}

\setlist[enumerate]{font = \normalfont}

\newcommand {\N}	{\mathrm{N}}

\newcommand {\R}	{\mathbb{R}}
\newcommand {\C}	{\mathbb{C}}
\newcommand {\E}	{\mathbb{E}}
\newcommand {\F}	{\mathbb{F}}

\newcommand {\T}	{\mathbb{T}}

\renewcommand{\d}{\, \mathrm{d}}

\DeclareMathOperator{\divH}{div_{\H}}

\newcommand{\sD}{\mathcal{D}}

\renewcommand{\Re}{\mathrm{Re}}
\renewcommand{\Im}{\mathrm{Im}}

\newcommand{\D}{\mathrm{D}}
\newcommand{\rN}{\mathrm{N}}
\renewcommand{\H}{\mathrm{H}}

\newcommand{\per}{\mathrm{per}}

\newcommand{\sigmabar}{\bar{\sigma}}

\newcommand{\atm}{\mathrm{a}}
\newcommand{\ocn}{\mathrm{o}}

\newcommand{\air}{\mathrm{a}}

\newcommand{\ssatm}[1]{{#1}^{\atm}}
\newcommand{\ssocn}[1]{{#1}^{\ocn}}

	\newcommand{\Omegaocn}{\ssocn{\Omega}}

	\newcommand{\Gau}{\Gamma_u^{\ocn}}
	\newcommand{\Gab}{\Gamma_b^{\ocn}}

	\newcommand{\dk}[1]{\partial_{#1}}
	\newcommand{\dt}{\dk{t}} 
	\newcommand{\dz}{\dk{z}} 
	
	\newcommand{\tr}{\mathrm{tr}}

	\newcommand{\eps}{\varepsilon}
	\renewcommand{\phi}{\varphi}
	\newcommand{\vatm}{\ssatm{v}}
	\newcommand{\vocn}{v^{\ocn}}
	\newcommand{\wocn}{w^{\ocn}}

	\newcommand{\uocn}{u^{\ocn}}

	\renewcommand{\bar}[1]{\overline{#1}}
	\newcommand{\vbar}{\bar{v}}
	\newcommand{\fbar}{\bar{f}}

	\renewcommand{\div}{\mathrm{div} \, }
	\newcommand{\nablaH}{\nabla_{\H}}
	\newcommand{\DeltaH}{\Delta_{\H}}

	\newcommand{\rC}{\mathrm{C}}
	\newcommand{\rL}{\mathrm{L}}
	\newcommand{\rW}{\mathrm{W}}
	\newcommand{\rH}{\H}
	\newcommand{\rB}{\mathrm{B}}

	\newcommand{\rLq}{\rL^q}
	\newcommand{\rLp}{\rL^p}

	\newcommand{\on}[1]{\text{ on } {#1}}

	\newcommand{\uair}{u^{\air}}
	\newcommand{\vair}{v^{\air}}

	\newcommand{\omegaair}{\omega^{\air}}

	\newcommand{\pair}{p^{\air}}

	\newcommand{\rhoocn}{\rho^{\ocn}}

	\newcommand{\Omegaair}{\Omega^{\air}}
	\newcommand{\rX}{\mathrm X}
    \newcommand{\Gammaairu}{\Gamma_u^{\atm}}
    \newcommand{\Gammaairb}{\Gamma_b^{\atm}}
    
    \newcommand{\wair}{w^\atm}

\title[The Incompressible-Compressible CAO problem]{The compressible-incompressible CAO problem introduced by Lions, Temam and Wang}

\author{Tarek Z\"{o}chling}
\address{Fachbereich Mathematik,
Technische Universit\"at Darmstadt,
Schlossgartenstr. 7, 64289 Darmstadt, Hessen, Germany}
\email{zoechling@mathematik.tu-darmstadt.de}
\thanks{Tarek Z\"ochling gratefully acknowledges the support by the Deutsche Forschungsgemeinschaft (DFG) through the Research Unit FOR~5528.}

\begin{document}

\begin{abstract}
\noindent
We study the (CAO) system of Lions-Temam-Wang, describing a compressible atmosphere coupled to an incompressible ocean through a fully nonlinear boundary condition representing wind stress at the interface.
In their foundational work, the system is reformulated in pressure coordinates under the assumption that the pressure is constant at the interface, leading to a purely incompressible system. In this article, we prove for the first time the strong well-posedness of the coupled compressible-incompressible system both locally and globally for small initial data.
\end{abstract}

\maketitle

\section{Introduction}
\label{sec:intro}

\noindent
The interaction between the atmosphere and the ocean has become a central focus in both the field of theoretical research and operational weather prediction. In this context, Lions, Temam, and Wang~\citep{LTW:93, LTW:95} developed a comprehensive model for atmosphere-ocean dynamics (CAO) that accounts for the effects of boundary layers near the interface. The model is based on a system of two coupled primitive equations. One describing atmospheric motion and the other governing oceanic dynamics, see also their earlier works~\citep{LTW:92a, LTW:92b}.
More specifically, the atmosphere is described by the \emph{compressible} primitive equations, which serve as a framework for large scale atmospheric behavior and the ocean is modeled using the \emph{incompressible} primitive equations, widely employed for understanding large-scale ocean currents. At the interface, the mechanical coupling process is realized via a fully nonlinear boundary condition that models the drag force exerted by the wind on the ocean surface. This interface condition plays a critical role in capturing the exchange of momentum between the two fluid domains and introduces significant mathematical challenges due to its nonlinearity. For a detailed description of the full model we refer to \autoref{sec: model}.

\medskip

\noindent
The primitive equations are widely considered a core model for representing geophysical flows, such as those in the ocean and atmosphere. They are derived from the Navier-Stokes or Boussinesq equations by simplifying the vertical component of motion under the assumption of hydrostatic balance. Since vertical velocity is difficult to observe directly, the hydrostatic approximation is commonly used. 

\medskip

\noindent
Readers interested in related material on geophysical fluid dynamics may consult the work of Chemin, Desjardins, Gallagher, and Grenier \citep{CDGG:06}. For comprehensive treatments of the Navier-Stokes equations, references include Galdi \citep{Gal:11}, Abe and Giga \citep{AG:13}, as well as Lemari\'e-Rieusset \citep{LR:16}.

\medskip

\noindent
In their pioneering work, Lions, Temam, and Wang established the existence of global, weak solutions to the (CAO)-system. They achieved this by formulating the compressible primitive equations in pressure coordinates (p-coordinates), showing that in this new coordinate system, the equations take the form of the classical primitive equations, with the incompressibility condition included. However,  the uniqueness of weak solutions and well-posedness of the (CAO)-system remained open problems. Recently, in \citep{BBHZ:25}, the well-posedness question was settled, showing that the (CAO)-system, subject to the pressure coordinate transform, is globally well-posed in the strong sense for arbitrarily large initial data in critical Besov spaces.

\medskip

\noindent
However, one limitation of the pressure coordinate transformation is the assumption of constant pressure at the interface, an assumption that lacks justification from a physical perspective, affects the accuracy of the wind force exerted by the atmosphere and most importantly simplifies the analysis by a margin, since the resulting system resembles the incompressible primitive equations.

\medskip

\noindent
In this article, we address this by studying the isothermal (CAO)-system without pressure coordinates. Specifically, we establish local, strong well-posedness for large initial data and global, strong well-posedness for small initial data. To the best of our knowledge, this is the first rigorous result concerning the isothermal (CAO)-system in its original, compressible-incompressible form. Our approach relies on the \emph{hydrostatic Lagrangian transform} introduced in \citep{HIRZ:25}, which differs from the classical Lagrangian transform, as it only considers the vertically averaged flow $\bar{v}^\air = \int_0^1 \vair(\cdot, \eta) \d \eta$. More precisely, we define the flow map $\mathrm{X}(t,y_\rH)$ by
\begin{equation}
\label{eq:flowpre}
\left \{
    \begin{aligned}
            \dt \mathrm{X}(t,y_\rH) &= \bar{v}^\air (t,\mathrm{X}(t,y_\rH)), \ t >0,\\
            \mathrm{X}(0,y_\rH) &= y_\rH,
    \end{aligned}
    \right. 
\end{equation}
where $y_\rH$ denotes the horizontal coordinate. This transformation is tailored to the anisotropic structure of the primitive equations and respects the replacement of the vertical momentum equation by the hydrostatic balance.
After the transformation, the linearized problem can be studied using optimal data results \citep{BBHZ:25, DHP:03} to control the nonlinear boundary conditions. Moreover, the global argument and the associated estimates are closely related in spirit to the classical theory for density-dependent (inhomogeneous) incompressible Navier--Stokes systems near a constant reference density, see, for instance, \citep{LS:78,AGZ:11,DM:12}.

\medskip

\noindent
Concerning the uncoupled incompressible primitive equations subject to standard Dirichlet and Neumann boundary conditions, global well-posed in 3D for arbitrary large $\rH^1$-data, was shown in the celebrated article \citep{CT:07} by Cao and Titi. In contrast, the compressible case is less understood. Under Neumann boundary conditions, both local and global strong well-posedness for small data without gravity were established in \citep{LT:20, LT:21}. In \citep{HIRZ:25}, these results are extended to incorporate the isothermal ideal gas law and gravity. Recently in \citep{HJMZ:25}, we also established the global well-posedness of the two-dimensional compressible primitive equations for arbitrary large data in $\rH^1$.
For further related results, see \citep{KZ:07, HK:16, GGHHK:20b}.

\medskip

\noindent
The structure of the article is as follows: In \autoref{sec: model} we present the coupled atmosphere ocean system as introduced in \citep{LTW:93}. In \autoref{sec:prelim} and \autoref{subsec: main} we state after some preliminaries the main results on local and global, strong well-posedness. \autoref{subsec:lagrangecor} recalls the hydrostatic Lagrangian transform and its basic properties from \citep{HIRZ:25}. Finally in \autoref{sec:lintheo} and \autoref{sec:global} the local and global well-posedness results are proven.

\section{The coupled atmosphere-ocean (CAO)-system}\label{sec: model}

\noindent
The coupled atmosphere-ocean (CAO) model describes a two-phase system where the atmosphere is represented by a compressible primitive equation, while the ocean is modeled by an incompressible primitive equation. These two subsystems are linked through a fully nonlinear boundary condition at the interface.
In the isothermal case, with all thermodynamic constants normalized to one, the coupled system is expressed as follows

\begin{equation}
	\left\{
	\begin{aligned}
        \dt \varrho^{\mathrm{a}} +\div  (\varrho^{\mathrm{a}} \uair)  &=0, &&\text{ in }\Omegaair \times (0,T), \\
		\varrho^{\mathrm{a}}\bigl (\dt  \vair + \uair \cdot \nabla \vair \bigr )  - \mu \Delta \vair- \mu' \nablaH \divH \vair +  \nablaH \pair &=0, &&\text{ in } \Omegaair \times (0,T),  \\
		  \dz \pair  &= -g\varrho^\air  , &&\text{ in } \Omegaair \times (0,T), \\ 
        \pair &=  \varrho^{\mathrm{a}}, &&\text{ in } \Omegaair \times (0,T),  \\
        \dt  \vocn - \Delta \vocn + \uocn \cdot \nabla \vocn + \nablaH \pi &=0, &&\text{ in } \Omegaocn \times (0,T),  \\
		  \dz \pi  &= -g, &&\text{ in } \Omegaocn \times (0,T), \\ 
		\div \uocn &= 0, &&\text{ in } \Omegaocn \times (0,T),\\
             \varrho^\air(0)=\varrho_0, \ \vair(0)=\vair_0,   \ \vocn(0) &= \vocn_0,
	\end{aligned}
	\right. 
	\label{eq:CAO}\tag{\textcolor{blue}{CAO}}
\end{equation}
Here $\Omegaair$ and $\Omegaocn$ are given by the cylindrical domains $\Omegaair = \T^2 \times (0,1)$ and $\Omegaocn = \T^2 \times (-1,0)$. 
The unknowns are the velocity of the atmosphere, denoted by $\uair = (\vair, \wair) : \Omegaair \rightarrow \R^3$, the atmospheric pressure $\pair : \Omegaair \rightarrow \R$ and the atmospheric density $\varrho^{\mathrm{a}} : \Omegaair \rightarrow \R$.
Concerning the oceanic part of the equations the unknowns are the velocity of the ocean~$\uocn = (\vocn , \wocn) : \Omegaocn \rightarrow \R^3$ as well as the pressure $\pi : \Omegaocn \rightarrow \R$.  Here we use $\nablaH$, $\divH$ for the horizontal gradient and the horizontal divergence, that is $\nablaH:=(\partial_x,\partial_y)^T$ and
$\divH:=\nablaH \cdot{}$. The Lam\'e coefficients $\mu$ and $\mu'$ are assumed to satisfy the standard conditions $\mu>0$ and  $\mu+ \mu'>0$. 



\medskip

\noindent
We now outline the mechanical coupling process used in the \eqref{eq:CAO}-model, following the description above. For a detailed discussion, we refer to \citep{LTW:95}.
The physical interaction between atmosphere and ocean is mediated at the air-sea interface \(\Gamma_i\). The boundary law at \(\Gamma_i\) governs the transfer of horizontal momentum and is obtained by balancing tangential (shear) stresses. The coupling between the atmosphere and ocean is derived from Newton's Laws of Motion. In particular, the primary driving mechanism for ocean currents is the wind drag. 

\medskip

\noindent
On the ocean side, the tangential shear stress reduces, at a flat interface, to the vertical shear of the horizontal oceanic velocity, which is given by
\[
(\text{\textcolor{blue}{oceanic shear stress}}) \;=\,\partial_z v^{\mathrm{o}}
\ \text{ on } \Gamma_i .
\]

On the atmospheric side, the tangential shear stress transmitted to the interface is proportional to the vertical shear of the horizontal atmospheric velocity and to the near-surface air density \(\varrho^{\mathrm{a}}\). With a positive constant \(c_1\) collecting the viscosity coefficients, we set
\[
(\text{\textcolor{blue}{atmospheric shear stress}}) \;=\; -\,c_1\,\varrho^{\mathrm{a}}\,\partial_z v^{\mathrm{a}}
\  \text{ on }  \Gamma_i .
\]

To represent the interfacial momentum input from the wind, we use the standard quadratic drag closure depending on the slip \(v^{\mathrm{a}}-v^{\mathrm{o}}\)
\begin{equation*}
  (\text{\textcolor{blue}{wind stress}}) \;= \;  c\,\varrho^{\mathrm{a}}(v^{\mathrm{a}}-v^{\mathrm{o}})\lvert v^{\mathrm{a}}-v^{\mathrm{o}}\rvert
    \  \text{ on } \Gamma_i,
\end{equation*}
with a drag coefficient \(c>0\). Continuity of tangential momentum flux across \(\Gamma_i\) identifies these three expressions and yields the fully nonlinear wind-stress coupling
\begin{equation}
     \label{eq:windcoupling}
-\,c_1\,\varrho^{\mathrm{a}}\,\partial_z v^{\mathrm{a}}
\;=\;
\partial_z v^{\mathrm{o}}
\;=\;
c\,\varrho^{\mathrm{a}}\,(v^{\mathrm{a}}-v^{\mathrm{o}})\,\lvert v^{\mathrm{a}}-v^{\mathrm{o}}\rvert
\ \text{ on } \Gamma_i .
\end{equation}

In isothermal regimes with constant pressure along the interface, the ideal-gas relation implies that \(\varrho^{\mathrm{a}}\) is constant on \(\Gamma_i\), which simplifies prefactors but leaves the mechanics unchanged. The interfacial tangential shear stress is continuous and is determined by a quadratic drag in the atmospheric-oceanic velocity difference. At the top of the atmosphere and at the bottom of the ocean, the boundaries are assumed to be impermeable and free-slip, meaning that no fluid crosses them and no tangential shear stress is imposed. Accordingly, the boundary conditions are
\begin{equation*}
        w^\air |_{z=1}= \partial_z v^\air |_{z=1}=0
\ \text{ and } \ w^\ocn |_{z=-1} =
\partial_z v^\ocn |_{z=-1}=0
\end{equation*}
\medskip


\noindent
The boundary conditions at the air-sea interface are central to any coupled model. While the complex formulation involving relative velocity is physically accurate, it is often simplified in theoretical studies. 
A common simplification of the wind stress $\tau = c \, \varrho^\air \, (\vair - \vocn) \, |\vair - \vocn|$
is to neglect the oceanic velocity, yielding 
\begin{equation*}
    \tau \approx \tau_0 = c\, \varrho^\air \, \vair\, |\vair|
    ,
\end{equation*}
see e.\ g., \citep{HellermanRosenstein:83}. This is justified by the scaling argument that ocean currents are typically much slower than winds, that is, $|\vocn| \ll |\vair|$.

\medskip

\noindent
This seemingly modest simplification neglects the so called \emph{eddy-killing effect}, removing a key feedback by which surface currents modulate the applied wind stress.
Because the stress then does work against the mesoscale currents, it extracts kinetic energy from eddies and damps their variance.
In practice, including the relative velocity in the stress term $\tau$ yields \emph{negative eddy wind work}, weakens the inverse cascade, and stabilizes western boundary currents, effects documented from early bulk formulations with current feedback, (Pacanowski~\citep{Pac:87}) through wind-work analyses of Scott-Xu~\citep{SX:09}, and in modern eddy-resolving studies (Renault et al.~\citep{Ren:19}).
Thus, while the simplification \(\tau\approx\tau_0\) may be convenient, it suppresses a key mechanical pathway by which the atmosphere regulates mesoscale ocean energetics.

\medskip

\noindent
For related results concerning linear boundary conditions, we refer e.\ g.\ to the work of Bresch and Simon  \citep{BS:01}. For related results about the Navier-Stokes equations with wind-driven boundary conditions, we 
refer to the work of Desjardins and Grenier \citep{DG:00}, Bresch and Simon \citep{BS:01}, Bresch, Guill\'en-Gonz\'alez, Masmoudi and Rodr\'iguez-Bellido \citep{BGMR:03} and Dalibard and Saint-Raymond \citep{DS:09}.

\section{Preliminaries}\label{sec:prelim}

\noindent
Let $T>0$ be a given time. We study the system \eqref{eq:CAO}, coupled through the condition \eqref{eq:windcoupling}, and impose the following boundary conditions on the top of the atmosphere, the bottom of the ocean, and along the lateral boundaries:

\begin{equation*} 
\begin{aligned} 
  \dz \vatm|_{z=1} = 0 \ &\text{ and } \ \wair|_{z=0,1} = 0 , \\ 
    \dz \vocn|_{z=-1} = 0  \ &\text{ and } \  \wocn|_{z=-1,0} = 0 ,\\
    \end{aligned} 
    \label{eq:bc1}
\end{equation*}

We deduce from \eqref{eq:CAO}$_3$ and \eqref{eq:CAO}$_4$ that 
\begin{equation}\label{eq:density}
    p^\air(t,x,y,z) = \varrho^\air(t,x,y,z) = \xi(t,x,y) \mathrm{e}^{-z},
\end{equation}
where we denote by $\xi$ the pressure evaluated at the surface, that is $$\xi(t,x,y) = \varrho^\air(t,x,y,z=0).$$ 
Moreover, we have set $g=1$, for simplicity of the notation. Our aim in the following is to exploit the explicit vertical dependence of the atmospheric density given by \eqref{eq:density}. This allows us to reduce the continuity equation to an effective two-dimensional equation for the surface pressure $\xi$ and, at the same time, to derive a representation formula for the vertical velocity $w^\air$.
To this end, we introduce in the atmospheric domain the change of variables
\begin{equation}\label{eq:ztransform}
    \mathbf z = \frac{1-\mathrm{e}^{-z}}{\delta}
    \iff
    z = \log\Bigl(\frac{1}{1-\delta \mathbf z}\Bigr),
    \qquad
    \delta = 1-\mathrm{e}^{-1} >0.
\end{equation}
We then define the transformed atmospheric unknowns by
\begin{equation*}
    \mathbf v^\air(t,x,y,\mathbf z)
    := v^\air\bigl(t,x,y,z(\mathbf z)\bigr),
    \qquad
    \mathbf w^\air(t,x,y,\mathbf z)
    := w^\air\bigl(t,x,y,z(\mathbf z)\bigr).
\end{equation*}
Since the transformation maps $(0,1)$ onto $(0,1)$, the domain $\Omega^\air$ remains unchanged. Moreover, by \eqref{eq:density}, the atmospheric pressure and density take the form
\begin{equation*}
    p^\air(t,x,y,\mathbf z)
    =
    \varrho^\air(t,x,y,\mathbf z)
    =
    \xi(t,x,y)\,(1-\delta \mathbf z).
\end{equation*}
Using \eqref{eq:density}, the atmospheric continuity equation can be written as
\begin{equation*}
    \mathrm e^{-z}\bigl(\partial_t \xi + \divH(\xi v^\air)\bigr)
    + \partial_z\bigl(\xi \mathrm e^{-z} w^\air\bigr)=0,
\end{equation*}
which, expressed in the new vertical variable $\mathbf z$, becomes
\begin{equation*}
    \partial_t \xi
    + \divH\bigl(\xi \mathbf v^\air\bigr)
    + \frac{1}{\delta}\,\partial_{\mathbf z}\bigl(\xi (1-\delta \mathbf z)\mathbf w^\air\bigr)
    =0.
\end{equation*}
Integrating with respect to $\mathbf z \in (0,1)$ and using the boundary conditions, we obtain
\begin{equation*}
    \partial_t \xi
    + \overline{\mathbf v}^{\air} \cdot \nablaH \xi
    + \xi \,\divH \overline{\mathbf v}^{\air}
    =0,
\end{equation*}
where
\[
\overline{\mathbf v}^{\air}
:=
\int_0^1 \mathbf v^\air(\cdot,\mathbf z)\d \mathbf z.
\]
Combining this with the transformed atmospheric momentum equation and the vertically averaged oceanic divergence condition, we arrive at the transformed system
\begin{equation}
	\left\{
	\begin{aligned}
        \partial_t \xi+ \overline{\mathbf v}^{\air} \cdot \nablaH\xi+ \xi\divH  \overline{\mathbf v}^{\air}&=0, &&\text{ in }\T^2 \times (0,T), \\
		\xi \bigl ( \partial_t \mathbf v^\air + \overline{\mathbf v}^{\air} \cdot  \nablaH \mathbf v^\air \bigr ) -\frac{\mu\DeltaH \mathbf v^\air}{1-\delta \mathbf z} -\mu \partial_{\mathbf z} \bigl ( \frac{1-\delta \mathbf z}{\delta^2} \partial_{\mathbf z} \mathbf v^\air \bigr ) -\frac{ \mu' \nablaH \divH \mathbf v^\air }{1-\delta \mathbf z}&= \mathbf F, &&\text{ in } \Omega^\air \times (0,T), \\
        \partial_t  v^\ocn - \Delta v^\ocn + u^\ocn \cdot \nabla v^\ocn + \nablaH \pi_s &= 0, &&\text{ in } \Omega^\ocn \times (0,T),  \\
		\div  \bar v^{\ocn} &= 0, &&\text{ in } \T^2 \times (0,T),\\
          \xi(0)=\xi_0 := \varrho_0(z=0), \ \mathbf v^\air(0)=v^\air_0, \ v^\ocn(0) &= v^\ocn_0,
	\end{aligned}
	\right.
	\label{eq:CAO1}
\end{equation}
where
\begin{equation*}
     \mathbf F(\xi,\mathbf v^\air)
     =
     - \xi \Bigl ((\mathbf v^\air-\overline{\mathbf v}^{\air}) \cdot \nablaH \mathbf v^\air
     + \frac{1-\delta \mathbf z}{\delta}\mathbf w^\air \partial_{\mathbf z}\mathbf v^\air \Bigr )
     - \nablaH \xi.
\end{equation*}
From now on, in order to simplify the notation, we drop the boldface and denote the transformed variables again by $z$, $v^\air$, and $w^\air$.

Here, the pressure $\pi$ is replaced by the surface pressure $\pi_s(t,x,y):= \pi(t,x,y,0)$ as it follows from \eqref{eq:CAO}$_6$ that $\pi = \pi_s -z$. Next, note that the vertical velocities $\wair$ and $\wocn$ are determined by their respective horizontal velocities and the density. More precisely, comparing \eqref{eq:CAO}$_1$ and \eqref{eq:CAO1}$_1$ as well as rearranging the divergence free condition \eqref{eq:CAO}$_7$ yields
\begin{equation}
    \label{eq:w}
    \xi w^\air = - \frac{\delta}{1- \delta z} \int_0^z \divH (\xi\vair-\xi \bar{v}^\air) \d z', \ z \in (0,1) \ \text{ and }\ \wocn = -\int_{-1}^z \divH \vocn \d z', \ z \in (-1,0),
\end{equation}
see also \citep{BBHZ:25, HIRZ:25}. For notational convenience, we now absorb the factor \(\frac{1-\delta z}{\delta}\) into the atmospheric vertical velocity and name it as before.
In particular, these representation formulas are consistent with the boundary conditions on the vertical velocities, since they immediately yield
\[
w^\air|_{z=0}=0
\qquad\text{and}\qquad
w^\ocn|_{z=-1}=0.
\]
Moreover, using \eqref{eq:CAO1}$_1$ and \eqref{eq:CAO1}$_4$, one also recovers the boundary conditions at the interface and the upper boundary, that is,
\[
w^\air|_{z=1}=0
\qquad\text{and}\qquad
w^\ocn|_{z=0}=0.
\]

The system is supplemented by the transformed boundary conditions
\begin{equation} 
\begin{aligned}\label{eq:bcCAO}
     \dz \vatm|_{z=1} = 0 \ &\text{ and } \ (\dz  \vatm)|_{z=0} = -\delta  (\vair-\vocn)|\vair-\vocn| , \\ 
   \dz \vocn|_{z=-1} = 0 \ &\text{ and } \ (\partial_z \vocn)|_{z=0} = \xi^{\air}(\vair-\vocn)|\vair-\vocn|. \\
    \end{aligned} 
\end{equation}

\begin{rem}\label{rem:systemeq}{\rm 
Let us note that the transformed system \eqref{eq:CAO1} is equivalent to the original system \eqref{eq:CAO} since the transform \eqref{eq:ztransform} is invertible for all times and the density can be recovered from \eqref{eq:density}.

}
\end{rem}

At this point, some words about our approach for solving \eqref{eq:CAO} and \eqref{eq:CAO1} are in order.
Let $s \in \R$ and $p,q \in (1,\infty)$. In the sequel, we
denote by $\rW^{s,q}(\Omega)$ the fractional Sobolev spaces, by $\rH^{s,q}(\Omega)$ the Bessel potential spaces, by $\rB^{s}_{qp}(\Omega)$ the Besov spaces, and by $\mathrm{F}_{pq}^s(\Omega)$ the Triebel--Lizorkin spaces, where $\Omega \subset \R^3$ is open. Moreover, the notation $\dot{\rH}^{s,q}(\Omega)$ refers to the corresponding homogeneous Sobolev space.
As usual, we set $\rH^{0,q}(\Omega) := \rL^q(\Omega)$ and note that $\rH^{s,q}(\Omega)$ coincides with the 
classical Sobolev space $\rW^{m,q}(\Omega)$ provided $s=m \in \N$. For more information on function spaces we refer to the monographs \citep{Tri:78, BCD:11, Ama:19}.

Next, we define \emph{hydrostatically solenoidal vector fields} by
\begin{equation*}
    \rL^q_{\sigmabar}(\Omega) = \overline{\{v \in \rC^\infty(\overline{\Omega};\R^2) \colon \divH \vbar = 0 \}}^{\| \cdot \|_{\rL^q(\Omega)}}.
\end{equation*}
Its role for the primitive equations parallels the one of the solenoidal vector fields for the Navier-Stokes equations. 
\section{Main results}\label{subsec: main}\mbox{} \\
\noindent
We now in the position to formulate our main theorems concerning the local well-posedness as well as the global well-posedness for small data of the \eqref{eq:CAO} system subject to the fully nonlinear boundary conditions \eqref{eq:bcCAO}. Our requirements on the set of initial data are collected in the following assumption.

\begin{asu}\label{assu:data0}
Let $p,q \in (2,\infty)$ such that $1 \neq \frac{1}{p} + \frac{1}{2q}$ and  $\frac{1}{2} \neq \frac{1}{p} + \frac{1}{2q}$ and assume that 
 \begin{enumerate}[(i)]
        \item $\xi_0 \in \rH^{1,q}(\T^2),\ \text{such that}\ M_1 \leq \xi_0 \leq M_2,\ \text{ for all }\ (x,y) \in \T^2\ \text{ and constants }\ M_1, M_2 >0$,
        \item $(v^\atm_0,v^\ocn_0) \in \rB_{qp}^{2(1-\nicefrac{1}{p})}(\Omega^\air;\R^2) \times \rB_{qp}^{2(1-\nicefrac{1}{p})}(\Omega^\ocn;\R^2) \cap \rLq_{\sigmabar}(\Omegaocn;\R^2)\ $ satisfying $(\dz v^\atm_0, \dz v^\ocn_0)|_{\Gammaairu \times \Gab}=0 \ $ and $\ (\dz v^\atm_0, \dz v^\ocn_0)|_{\Gammaairb \times \Gau}=(\delta  (\vocn_0 - \vair_0) |\vair_0 - \vocn_0|,\xi^{\air}_0  (\vair_0 - \vocn_0) |\vair_0 - \vocn_0|) \ $ if $\ \frac{1}{2q}+ \frac{1}{p} < \frac{1}{2}$.
    \end{enumerate}
\end{asu}

\begin{thm}[local strong well-posedness]\label{thm:localwpCAO}\mbox{}\\
   Let $T>0$ and assume that $(\xi_0, \vocn_0,\vatm_0)$ satisfy \autoref{assu:data0}.
    Then there is $0 < a =a(\vair_0, \vocn_0) \leq T$ such that the system \eqref{eq:CAO} admits a unique, strong solution $(\xi^{\air}, \vair, \vocn,\pi_s)$ subject to the boundary conditions \eqref{eq:bcCAO} satisfying
    \begin{equation}
        \begin{aligned}
            \label{eq:solCAO}
            \xi&\in \rH^{1,p}(0,a;\rL^{q}(\T^2)) \cap  \rL^\infty(0,a;\rH^{1,q}(\T^2)), \\
            \vair &\in \rH^{1,p}(0,a; \rLq(\Omegaair,\R^2)) \cap \rLp(0,a;\rH^{2,q}(\Omegaair;\R^2)),  \\
            \vocn &\in  \rH^{1,p}(0,a; \rLq(\Omegaocn;\R^2)) \cap \rLp(0,a;\rH^{2,q}(\Omegaocn;\R^2)), \\
            \pi_s &\in \rLp(0,a;\dot \rH^{1,q}(\T^2))
        \end{aligned}
    \end{equation}
    and there are constants $M_1^\ast$, $M_2^\ast>0$ such that $M_1^\ast \leq \xi \leq M_2^\ast$ for all $t \in (0,a)$, $(x,y) \in \T^2$.
\end{thm}
Next we state the result concerning global well-posedness for small initial data.
\begin{thm}[global strong well-posedness for small data]
    \label{thm:globalWP}\mbox{}\\
Let $T>0$. Assume that $(\xi_0, \vatm_0.\vocn_0)$ satisfy \autoref{assu:data0} and in addition 
    \begin{equation*}
        \begin{aligned}
           \| \xi_0-\bar{\xi} \|_{\rH^{1,q}(\T^2)} + \| v^\air_0 \|_{\rB^{2(1-\nicefrac{1}{p})}_{qp}(\Omegaair)} +   \| v^\ocn_0 \|_{\rB^{2(1-\nicefrac{1}{p})}_{qp}(\Omegaocn)} \leq \delta,
        \end{aligned}
    \end{equation*}
    for some $\delta>0$ and $\bar{\xi} >0$. Then the system \eqref{eq:CAO} admits a unique, strong solution $(\xi,\vair, \vocn, \pi_s)$ subject to the boundary conditions \eqref{eq:bcCAO} satisfying \eqref{eq:solCAO} with $a=T$ as well as 
   \begin{equation*}
        \begin{aligned}
            &\| \xi-\bar{\xi} \|_{\rL^\infty(0,T;\rH^{1,q}(\T^2))} + \| \nablaH \xi \|_{\rL^{\infty}(0,T; \rLq(\T^2))} + \|  \dt \xi \|_{\rLp(0,T; \rLq(\T^2))}\\& +\|  \vair \|_{\rLp(0,T;\rH^{2,q}(\Omegaair))} + \|  \dt \vair\|_{\rLp(0,T ; \rLq(\Omegaair))} + \|  \vair \|_{\rL^\infty(0,T; \rB_{qp}^{2(1-\nicefrac{1}{p})}(\Omegaair))} \\
            &+ \|  \vocn \|_{\rLp(0,T;\rH^{2,q}(\Omegaocn))} + \|  \dt \vocn\|_{\rLp(0,T ; \rLq(\Omegaocn))} + \| \vocn \|_{\rL^\infty(0,T;\rB_{qp}^{2(1-\nicefrac{1}{p})}(\Omegaocn))} \\
           &+ \|  \pi_s \|_{\rLp(0,T; \rH^{1,q}(\T^2)) }\leq C \delta.
        \end{aligned}
      \end{equation*}
    Moreover, it holds that $\xi(t,x,y) \geq \nicefrac{\bar{\xi}}{2}$ for all $t \in [0,T]$, $(x,y) \in \T^2$.
\end{thm}

\begin{rem}
    \label{rem:pressure}{\rm
  \autoref{thm:localwpCAO} also holds true when considering a power law connecting the pressure and the density of the atmosphere, i.\ e., if \eqref{eq:CAO}$_4$ is replaced by the equation $\pair=(\varrho^\air)^\gamma$, for $\gamma>1$. Moreover, in case of $g=0$, both \autoref{thm:localwpCAO} and \autoref{thm:globalWP} hold true when considering a generalized pressure law, meaning that 
        \begin{equation*}
    \pair = P(\varrho^\air), \ \text{ where }P \in \rC^\infty \ \text{ with }\ c_1 \leq  P'(s) \leq c_2, \ \forall s \in \Im \ \varrho^\air,
\end{equation*}
for given positive constants $c_1$ and $c_2$, see also \citep{HIRZ:25}.
        
}
\end{rem}

\section{Hydrostatic Lagrangian coordinates}\label{subsec:lagrangecor}

\noindent
In this section we recall the hydrostatic Lagrangian transformation introduced in \citep{HIRZ:25}. Let $T>0$ fixed and consider the characteristic along $\bar{v}^\air$, i.\ e., the solution $\mathrm{X}$ of the differential equation given by
\begin{equation}
\label{eq:flow}
\left \{
    \begin{aligned}
            \dt \mathrm{X}(t,y_\H) &= \bar{v}^\air (t,\mathrm{X}(t,y_\H)), \quad t \in (0,T) \\
            \mathrm{X}(0,y_\H) &= y_\H,
    \end{aligned}
    \right. 
\end{equation}
where $y_\H = (y_1,y_2) \in \T^2$. The mapping $\mathrm{X}(t, \cdot) : \T^2 \rightarrow \T^2$ defines a well-defined $\rC^1$-diffeomorphism, assuming sufficient regularity of $\bar{v}^\air$. Solving \eqref{eq:flow} leads to
\begin{equation}\label{eq:diffeo X Euler Lagrange sea ice para-hyper}
    \rX(t,y_\H) = y_\H + \int_0^t \ \bar{v}^\air(s,\rX(s,y_\H)) \d s, 
\end{equation}
for $y_\H \in \T^2$ and $t \in (0,T)$. The motivation for using this specific transformation arises from the anisotropic structure of the solution $\uair =(\vair, \wair)$. In particular, as seen from \eqref{eq:w}, the vertical velocity $w^\air$ behaves horizontally like $\divH(\xi \vair)$ and is thus less regular. 
Denoting by $\mathrm{Y}(t,\cdot)$, the inverse of $\rX(t,\cdot)$, by $\mathrm{Z} = [\nablaH \rX]^{-1}$ and by $\mathrm J = \det \nablaH \rX$. 
We recall the following properties of the hydrostatic Lagrangian transformation from \citep{HIRZ:25}.

\begin{lem}\label{lem:ests of trafo}
Let $T>0$ and $p,q \in (2,\infty)$.
Assume that
\begin{equation*}
    \begin{aligned}
        v^\air \in \E^\atm_{1,T} := \rH^{1,p}(0,T;\rLq(\Omegaair;\R^2)) \cap \rLp (0,T;\rH^{2,q}(\Omegaair;\R^2)) \ \text{ with } \  \| v^\air \|_{\E_{1,T}^\atm} \leq R,
    \end{aligned}
\end{equation*}
for a constant $R>0$ independent of $T$.
Let $\rX$, $\mathrm{Z}$ and $\mathrm{J}$ be the flow quantities associated with the averaged velocity $\bar{v}^\air$. 
Then the following statements hold true.

\begin{enumerate}[(a)]
    \item There is $C > 0$, independent of $T$, such that
    \begin{equation*}
        \sup_{t \in (0,T)} \|\nablaH \rX(t,\cdot) - \mathbb{I}_2 \|_{\rL^\infty(\T^2)}
        \le C \cdot \sup_{t \in (0,T)} \|\nablaH \rX(t,\cdot) - \mathbb{I}_2 \|_{\rH^{1,q}(\T^2)}
        \le C T^{\nicefrac{1}{p'}}.
    \end{equation*}
    In particular, there is $T'$ such that for $0 < t < T \le T'$, we have
    \[
    \|\nablaH \rX(t,\cdot) - \mathbb{I}_2 \|_{\rL^\infty(\T^2)} \le \nicefrac{1}{2},
    \]
    so $\nablaH \rX(t,\cdot)$ is invertible, and $\mathrm{Z} = [\nablaH \rX]^{-1}$ is thus well-defined.

    \item For $0 < T \le T'$, with $T'$ as made precise in~(a), there is $C > 0$, independent of $T$, so that
    \begin{equation*}
        \begin{aligned}
            \| \mathrm{Z} \|_{\rH^{1,p}(0,T;\rH^{1,q}(\T^2))}
            + \| \mathrm{Z} \|_{\rL^\infty(0,T;\rH^{1,q}(\T^2))}
            + \|\nablaH \rX \|_{\rH^{1,p}(0,T;\rH^{1,q}(\T^2))}
            + \|\nablaH \rX \|_{\rL^\infty(0,T;\rH^{1,q}(\T^2))}
            \le C.
        \end{aligned}
    \end{equation*}

    \item Similar estimates as in (a) are also valid for $\mathrm{Z}$ as well as $\mathrm{Z}^\top$, i.\ e.,
    \begin{equation*}
        \sup_{t \in (0,T)} \| \mathrm{Z}(t,\cdot) - \mathbb{I}_2\|_{\rH^{1,q}(\T^2)}
        + \sup_{t \in (0,T)} \| \mathrm{Z}^\top(t,\cdot) - \mathbb{I}_2 \|_{\rH^{1,q}(\T^2)}
        \le C T^{\nicefrac{1}{p'}}
    \end{equation*}
    holds for a constant $C>0$ independent of $T$.

    \item For $0<T\le T'$, it holds that $\mathrm{J}(t,y_\rH) \geq c >0$ for all $t \in (0,T)$ and $y_\rH \in \T^2$. Moreover, there is $C>0$, independent of $T$, such that
    \begin{equation*}
        \| \mathrm J \|_{\rH^{1,p}(0,T;\rH^{1,q}(\T^2))}
        + \| \mathrm J \|_{\rL^\infty(0,T;\rH^{1,q}(\T^2))}
        \le C \ \text{ and } \  \sup_{t\in(0,T)} \| \mathrm J(t,\cdot)-1 \|_{\rH^{1,q}(\T^2)}
        \le C T^{\nicefrac{1}{p'}}.
    \end{equation*}

    \item For all $j,k,l \in \{1,2\}$, there is $C > 0$, independent of $T$, such that
    \begin{equation*}
        \Bigl\| \frac{\partial \mathrm{Z}_{l,j}}{\partial y_k} \Bigr\|_{\rL^\infty(0,T;\rL^q(\T^2))} \le C T^{\nicefrac{1}{p'}}.
    \end{equation*}
\end{enumerate}
\end{lem}

In the following, we express the oceanic system in the same transformed notation as the atmospheric part. Although, strictly speaking, this is not necessary, one may also keep the oceanic equations in the original variables and transform only the coupling terms on the interface, adopting a common notation considerably simplifies the presentation and avoids repeated switching between coordinate systems.
Considering the transformation
\[
(t,y_\H,z)\mapsto (t,\rX(t,y_\H),z)
\]
as well as the new quantities
\begin{equation*}
    \begin{aligned}
            \zeta(t,y_\H) := \xi(t, \mathrm{X}(t,y_\H)), \quad V^\air(t,y_\H,z) &:= \vair (t, \mathrm{X}(t,y_\H),z), \quad   V^\ocn(t,y_H,z) := \vocn (t, \mathrm{X}(t,y_\H),z)\ \text{ and } \\  \Pi(t,y_\H) &:= \pi_s(t,\rX(t,y_\H)).
    \end{aligned}
\end{equation*}
Then the transformed surface density satisfies
\[
\partial_t \zeta
+
\zeta\,\tr \, (\mathrm Z^\top \nablaH \bar V^\air)=0,
\]
whereas the Jacobian
solves
\[
\partial_t \mathrm J
=
\mathrm J\,\tr\!\bigl(\mathrm Z^\top \nablaH \bar V^\air\bigr) \ \text{ with }\
\mathrm J(0,\cdot)=1.
\]
Hence
\[
\partial_t(\zeta \mathrm J)=0
\ \text{ and therefore }\
\zeta(t,y_\H)\,\mathrm J(t,y_\H)=\xi_0(y_\H).
\]
In particular, the transformed density is completely determined by the initial data and the Jacobian of the flow, namely
\begin{equation}\label{eq:zetaJ}
    \zeta(t,y_\H)=\frac{\xi_0(y_\H)}{\mathrm J(t,y_\H)}.
\end{equation}
The remaining system is transformed to
\begin{equation}
	\left\{
	\begin{aligned}
		\partial_t  V^\air
        - \frac{\mu \DeltaH V^\air }{(1-\delta z)\xi_0}
        - \frac{\mu}{\xi_0} \partial_z \bigl ( \frac{1-\delta z}{\delta^2} \partial_z V^\air \bigr )
        -  \frac{\mu' \nablaH \divH V^\air}{(1-\delta z)\xi_0}
        &= F_1(\zeta, V^\air), &&\text{ in } \Omega^\air \times (0,T),  \\
        \partial_t  V^\ocn - \Delta V^\ocn + \nablaH \Pi
        &= F_2(\Pi, V^\ocn), &&\text{ in } \Omega^\ocn \times (0,T),  \\
		\divH  \bar{V}^{\ocn}
        &= F_3(V^\ocn), &&\text{ in } \T^2 \times (0,T),\\
        V^\air(0)&=v^\air_0, \quad V^\ocn(0)=v^\ocn_0,
	\end{aligned}
	\right.
	\label{eq:CAOatm}
\end{equation}
where $\zeta$ is given by \eqref{eq:zetaJ} and the right-hand sides are given by
\begin{equation*}
\begin{aligned}
(F_1)_i
&=
\frac{\mu}{(1-\delta z)\xi_0}
\Bigl(
\sum_{l,j,k}
\frac{\partial^2 V^\air_i}{\partial y_l \partial y_k}
\bigl(\mathrm Z_{k,j}-\delta_{k,j}\bigr)\mathrm Z_{l,j}
+
\sum_{l,k}
\frac{\partial^2 V^\air_i}{\partial y_l \partial y_k}
\bigl(\mathrm Z_{l,k}-\delta_{l,k}\bigr)
+
\sum_{l,j,k}
\mathrm Z_{l,j}
\frac{\partial V^\air_i}{\partial y_k}
\frac{\partial \mathrm Z_{k,j}}{\partial y_l}
\Bigr) \\
&\quad
+
\frac{\mu'}{(1-\delta z)\xi_0}
\Bigl(
\sum_{l,j,k}
\frac{\partial^2 V^\air_j}{\partial y_l \partial y_k}
\bigl(\mathrm Z_{k,j}-\delta_{k,j}\bigr)\mathrm Z_{l,i}
+
\sum_{l,j}
\frac{\partial^2 V^\air_j}{\partial y_l \partial y_j}
\bigl(\mathrm Z_{l,i}-\delta_{l,i}\bigr)
+
\sum_{l,j,k}
\mathrm Z_{l,i}
\frac{\partial V^\air_j}{\partial y_k}
\frac{\partial \mathrm Z_{k,j}}{\partial y_l}
\Bigr) \\
&\quad
+
\Bigl(1-\frac{\zeta}{\xi_0}\Bigr)(\partial_t V^\air)_i
-\frac{\zeta}{\xi_0}
\Bigl(
\sum_{l,j}
(V^\air-\bar V^\air)_l
\frac{\partial V^\air_i}{\partial y_j}\,
\mathrm Z_{l,j}
+
\frac{1-\delta z}{\delta}\,
W^\air \frac{\partial V^\air_i}{\partial z}
\Bigr)
-\frac{(\mathrm Z^\top\nablaH \zeta)_i}{\xi_0},
\end{aligned}
\end{equation*}
as well as
\begin{equation*}
\begin{aligned}
(F_2)_i
&=
\sum_{l,j,k}
\frac{\partial^2 V^\ocn_i}{\partial y_l \partial y_k}
\bigl(\mathrm Z_{k,j}-\delta_{k,j}\bigr)\mathrm Z_{l,j}
+
\sum_{l,k}
\frac{\partial^2 V^\ocn_i}{\partial y_l \partial y_k}
\bigl(\mathrm Z_{l,k}-\delta_{l,k}\bigr)
+
\sum_{l,j,k}
\mathrm Z_{l,j}
\frac{\partial V^\ocn_i}{\partial y_k}
\frac{\partial \mathrm Z_{k,j}}{\partial y_l}
\\ &\quad
-
\sum_{l,j}
(V^\ocn-\bar V^\air)_l
\frac{\partial V^\ocn_i}{\partial y_j}\,
\mathrm Z_{l,j}
-
W^\ocn \frac{\partial V^\ocn_i}{\partial z}
-
\sum_j
\frac{\partial \Pi}{\partial y_j}
\bigl(\mathrm Z_{j,i}-\delta_{j,i}\bigr),
\end{aligned}
\end{equation*}
and
\begin{equation*}
    F_3
    =-
    \nablaH \bar{V}^\ocn : \left [\mathrm{Z}^\top  - \mathbb{I}_2 \right ]
    =
    -\divH \bigl ( (\mathrm{J}\mathrm{Z}^\top-\mathbb{I}_2 ) \bar{V}^\ocn \bigr ).
\end{equation*}

Finally, we remark that in hydrostatic Lagrangian variables the boundary conditions take the form of \eqref{eq:bcCAO}, with all quantities replaced by their transformed counterparts and with the density factor given by \eqref{eq:zetaJ}.

\section{Linear Theory and local-wellposedness}\label{sec:lintheo}

\noindent
Let $T>0$. To analyze the transformed system \eqref{eq:CAOatm}, we consider the following two decoupled linear problems. 
Concerning the velocity of the atmosphere $V^\air \colon \Omegaair \rightarrow \R^2$ consider the hydrostatic Lam\'{e} system
\begin{equation}
	\left\{
	\begin{aligned}
		\dt  V^\air - \frac{\mu \DeltaH V^\air }{(1-\delta z)\xi_0}- \frac{\mu}{\xi_0} \dz \bigl ( \frac{1-\delta z}{\delta^2} \dz V^\air \bigr ) -  \frac{\mu' \nablaH \divH V^\air}{(1-\delta z)\xi_0}  &=  f_1, &&\text{ in } \Omegaair \times (0,T),  \\
        \dz V^\air   &= 0  , &&\text{ on } \T^2 \times \{1\} \times (0,T), \\ 
		 \dz V^\air   &= b_1  , &&\text{ on } \T^2 \times \{0\} \times (0,T), \\ 
        V^\air(0)  &=v^\air_0 ,  &&\text{ in } \Omegaair,
	\end{aligned}
	\right. 
	\label{eq:Lame1}
\end{equation}
and for the velocity of the ocean $V^\ocn \colon \Omegaocn \rightarrow \R^2$ and the pressure of the ocean $\Pi \colon \T^2 \rightarrow \R$ consider the hydrostatic Stokes problem 
\begin{equation}
	\left\{
	\begin{aligned}
		\dt V^\ocn - \Delta V^\ocn+ \nablaH \Pi &=f_2, &&\text{ in } \Omegaocn \times (0,T), \\
		\dz \Pi  &= 0, &&\text{ in } \Omegaocn \times (0,T), \\
		\divH \bar{V}^\ocn &= f_3, &&\text{ in } \Omegaocn \times (0,T), \\
		\dz V^\ocn &=b_2, &&\text{ on } \T^2 \times \{0\} \times (0,T), \\
		\dz  V^\ocn &= 0, &&\text{ on } \T^2 \times \{-1\}\times (0,T), \\
        V^\ocn(0) &= v^\ocn_0, &&\text{ in } \Omegaocn. \\
	\end{aligned}
	\right. 
	\label{eq:hstokes1}
\end{equation}
The following assumptions are crucial for solving the boundary value problems \eqref{eq:Lame1} and \eqref{eq:hstokes1} subject to inhomogeneous data.

\begin{asu}\label{assu:data}
Let $p,q \in (2,\infty)$ such that $1 \neq \frac{1}{p} + \frac{1}{2q}$, $\frac{1}{2} \neq \frac{1}{p} + \frac{1}{2q}$ and assume that $f_1$, $f_2$, $f_3$, $b_1$, $b_2$, $\xi_0$, $v^\atm_0$, $v^\ocn_0$, satisfies 
\begin{enumerate}[(i)]
   \item $\xi_0 \in \rH^{1,q}(\T^2),\ \text{such that}\ M_1 \leq \xi_0 \leq M_2,\ \text{ for all }\ (x,y) \in \T^2\ \text{ and constants }\ M_1, M_2 >0$,
        \item $(v^\atm_0,v^\ocn_0) \in \rB_{qp}^{2(1-\nicefrac{1}{p})}(\Omega^\air;\R^2) \times \rB_{qp}^{2(1-\nicefrac{1}{p})}(\Omega^\ocn;\R^2) \ $ satisfying the compatibility conditions specified in \autoref{assu:data0} as well as $\divH \bar{v}_0^\ocn = f_4(0)$ in $\mathcal{D}'(\T^2)$, i.\ e., in the sense of distributions, 
    \item $(f_1,f_2) \in   \rLp(0,T;\rLq(\Omegaair;\R^2)) \times \rLp(0,T;\rLq(\Omegaocn;\R^2)) $,
    \item $f_3 \in \rH^{1,p}(0,T;\Dot{\rH}^{-1,q}(\T^2)) \cap \rLp(0,T;\rH^{1,q}(\T^2)) \eqqcolon \mathbb{G}_{1,T}$,
     \item $b_i \in \mathrm{F}_{pq}^{\nicefrac{1}{2}-\nicefrac{1}{2q}}(0,T;\rLq(\T^2;\R^2)) \cap \rLp(0,T;\rB_{qq}^{1 - \nicefrac{1}{q}}(\T^2;\R^2)) \eqqcolon \F_{\rN,T},\ $ for $i\in \{ 1,2\}$.
	\end{enumerate}
\end{asu}

The operator associated with \eqref{eq:Lame1} is the \emph{hydrostatic Lam\'e operator}, introduced and studied in \citep{HIRZ:25}. In particular, its symbol is parameter-elliptic with angle $0$, and the Lopatinski--Shapiro conditions are satisfied.
Hence, by \citep{DHP:07}, optimal data results concerning
\eqref{eq:Lame1} are known for the half-space
and can be adapted to the layer $\Omegaair_L = \R^2 \times (0,1)$
via the same localization procedure.
Finally, these results carry over to the domain
$\Omegaair = \T^2 \times (0,1)$
by extending periodic functions on $\T^2$ to the full space $\R^2$,
as demonstrated in \citep[p.~1082]{HK:16}.
Together with the considerations in \citep{BBHZ:25} concerning \eqref{eq:hstokes1} we obtain the following result.
\begin{lem}\label{lem:Lame equation}
 Let $T>0$ and $p,q$, $(f_1,f_2,f_3,b_1,b_2,\xi_0,\vair_0, \vocn_0)$ as in \autoref{assu:data}.

    $\mathrm{(a)}$ The problem \eqref{eq:Lame1} admits a unique, strong solution 
\begin{equation*}
    \begin{aligned}
        V^\air &\in \E^\air_{1,T} := \rH^{1,p}(0,T;\rLq(\Omegaair;\R^2)) \cap \rLp(0,T;\rH^{2,q}(\Omegaair;\R^2))
    \end{aligned}
\end{equation*}
satisfying
\begin{equation*}\label{eq: vair est}
    \| V^\air \|_{\E^\air_{1,T}} \le C\bigl(\| f_1\|_{ \rLp(0,T;\rLq(\Omegaair)) }  + \| \vair_0 \|_{\rB_{qp}^{2(1-\nicefrac{1}{p})}(\Omega^\air)} +  \| b_1 \|_{\F_{\N,T}}\bigr),
\end{equation*}
for some $C>0$.


    \noindent
    $\mathrm{(b)}$  The problem \eqref{eq:hstokes1} admits a unique, strong solution 
    \begin{equation*}
        V^\ocn \in \E^\ocn_{1,T}\coloneqq\rH^{1,p}(0,T;\rLq(\Omegaocn;\R^2)) \cap \rLp(0,T;\rH^{2,q}(\Omegaocn;\R^2)) \ \text{ and } \ \Pi_s \in \rLp(0,T; \dot {\rH}^{1,q}(\T^2) ) , 
    \end{equation*}
    satisfying
    \begin{equation*}\label{eq: est ocn}
    \| V^\ocn \|_{\E^\ocn_{1,T}} + \| \nablaH \Pi_s \|_{\rLp (0,T;\rLq(\T^2))} \le C\bigl(\| f_2\|_{\rLp(0,T;\rLq(\Omegaocn))}  +\| f_3 \|_{ \mathbb{G}_{1,T}}
    + \| \vocn_0 \|_{\rB_{qp}^{2(1-\nicefrac{1}{p})}(\Omega^\ocn) } + \| b_2 \|_{\F_{\N,T}}\bigr),
\end{equation*}
for some $C>0$.

\end{lem}

The proof of \autoref{thm:localwpCAO} is based on the contraction mapping principle, which in turn relies on \autoref{lem:Lame equation}, estimates for the nonlinear terms $F_1$, $F_2$, $F_3$, and the nonlinear coupling conditions given in \eqref{eq:bcCAO}. For simplicity of the notation,  we introduce the full data  and solution space by
\begin{equation*}
    \begin{aligned}
    \E_{0,T} &:=  \rLp(0,T;\rLq(\Omegaair;\R^2)) \times \rLp(0,T;\rLq(\Omegaocn;\R^2)) \times \mathbb{G}_{1,T}, \\
        \E_{1,T} &:=  \E^\atm_{1,T} \times \E^\ocn_{1,T} \times  \rLp(0,T; \dot {\rH}^{1,q}(\T^2) ) 
    \end{aligned}
\end{equation*}
and define the boundary operator 
\[
    B(V^\atm,V^\ocn) := \bigl( (V^\atm - V^\ocn)\,|V^\atm- V^\ocn| \bigr)|_{z=0},
\]
on the common flat interface.
The estimates of the nonlinear terms are given in the following lemma.
\begin{lem}
    \label{lem:nonlinearestimates}
     Let $T>0$ and $p,q \in (2,\infty)$. Assume that 
     \begin{equation*}
         (V^\atm, V^\ocn, \Pi) \in \E_{1,T} \ \text{ with } \ \| (V^\atm, V^\ocn, \Pi) \|_{ \E_{1,T}} \leq R,
     \end{equation*}
     for some constant $R>0$ independent of $T$. Then the following estimates hold true.
      \begin{equation*}
        \begin{aligned}
        \| (F_1(\zeta,V^\atm) , F_2(V^\ocn,\Pi) , F_3(V^\ocn)  ) \|_{\E_{0,T}}\leq C
( T^\delta + R^2 ) \ \text{ and } \  \| \zeta B(V^\atm,V^\ocn) \|_{\F_{\N,T}} \leq C R^2.
        \end{aligned}
        \end{equation*}
 for some $\delta >0$.
\end{lem}

\begin{proof}
By \eqref{eq:zetaJ}, we have
\[
\zeta=\frac{\xi_0}{\mathrm J}.
\]
Since $\mathrm J$ is bounded away from zero by Lemma~\ref{lem:ests of trafo}(d), it follows from the Banach algebra property of $\rH^{1,q}(\T^2)$ for $q>2$, the estimates for $\mathrm J$ and $\mathrm J-1$ from Lemma~\ref{lem:ests of trafo}(d), and standard composition results that
\begin{equation}\label{eq:zetaestimate}
    \|\zeta\|_{\rL^{\infty}(0,T;\rH^{1,q}(\T^2))}
    \le C, \quad 
    \|\zeta\|_{\rLp(0,T;\rH^{1,q}(\T^2))}
    \le CT^{\nicefrac{1}{p}},
    \ \text{ and } \ 
    \|\zeta-\xi_0\|_{\rL^\infty(0,T;\rH^{1,q}(\T^2))}
    \le C T^{\nicefrac{1}{p'}},
\end{equation}
where $C>0$ is independent of $T$. We will frequently use these bounds throughout the proof.
 We begin by estimating the first nonlinear term $F_1$ appearing on the right-hand side of the transformed system \eqref{eq:CAOatm}. Observe that
        \begin{equation*}
        \begin{aligned}
          &\big \|  \sum_{l,j,k} \frac{\partial^2 V^\atm_i}{\partial y_l \partial y_k} (\mathrm{Z}_{k,j} - \delta_{k,j}) \mathrm{Z}_{l,j} \big \|_{\rLp(0,T;\rLq(\Omegaair))}
          \\&\leq C \| V^\atm \|_{\rLp(0,T;\rH^{2,q}(\Omegaair))} \| \mathrm{Z} - \mathbb{I}_2 \|_{\rL^\infty (0,T;\rH^{1,q}(\T^2))} \| \mathrm{Z} \|_{\rL^\infty (0,T;\rH^{1,q}(\T^2))} \\
          &\leq C T^{\nicefrac{1}{p'}},
        \end{aligned}
     \end{equation*}   
      and similarly
\begin{equation*}
\begin{aligned}
          \big \|  \sum_{l,k}\frac{\partial^2 V_i}{\partial y_l \partial y_k}(\mathrm{Z}_{l,k}- \delta_{k,l}) \big \|_{\rLp(0,T;\rLq(\Omega))}
            &\leq C \| V^\atm \|_{\rLp(0,T;\rH^{2,q}(\Omega))} \| \mathrm{Z} - \mathbb{I}_2 \|_{\rL^\infty (0,T;\rH^{1,q}(\T^2))}\\&\leq C T^{\nicefrac{1}{p'}}.
        \end{aligned}
    \end{equation*}
    Furthermore, we obtain
    \begin{equation*}
        \begin{aligned}
          &\big \| \sum_{l,j,k} \mathrm{Z}_{l,k} \frac{\partial V^\atm_i}{\partial y_k} \frac{\partial \mathrm{Z}_{k,j}}{\partial y_l} \big \|_{\rLp(0,T;\rLq(\Omegaair))}
          \\&\leq C\| \mathrm{Z} \|_{\rL^\infty (0,T;\rL^{\infty}(\T^2))} \| \nablaH V^\atm \|_{\rLp(0,T;\rH^{1,q}(\Omegaair))}  \sum_{l,j,k}
          \big \| \frac{\partial \mathrm{Z}_{k,j}}{\partial y_l} \big \|_{\rL^\infty (0,T;\rLq(\T^2))} \\
            &\leq C T^{\nicefrac{1}{p'}}.
        \end{aligned}
    \end{equation*}
    Concerning the transformed time derivative, we find
      \begin{equation*}
\begin{aligned}
         \big\| \big (1 -   \frac{\zeta}{\xi_0} \big )(\dt V^\atm ) \big\|_{\rLp(0, T; \rLq(\Omegaair))} &\leq C\| \zeta - \xi_0 \|_{\rL^\infty(0,T;\rH^{1,q}(\T^2))} \| V^\atm\|_{\rH^{1,p}(0,T;\rLq(\Omegaair))}
        \\& \leq C T^{\nicefrac{1}{p'}}.
\end{aligned}
    \end{equation*}
    To estimate the nonlinear term, we first observe that for all $q>2$
    \begin{equation*}
        \begin{aligned}
            \| V^\atm \cdot \nablaH V^\atm \|_{\rLq(\Omegaair)} &\leq C\| V^\atm \|_{\rL^{3q}(\Omegaair)} \| V^\atm \|_{\rH^{1,\nicefrac{3q}{2}}(\Omegaair)} \le C \cdot \| V^\atm \|^2_{\rH^{1+ \nicefrac{1}{q},q}(\Omegaair)}.
        \end{aligned}
    \end{equation*}
    It follows that
    \begin{equation*}
        \begin{aligned}
              &\|\zeta \sum_{l,j}\Tilde{V}^\atm_l \frac{\partial V^\atm_i}{\partial y_j} \mathrm{Z}_{l,j}  \|_{\rLp(0,T;\rLq(\Omegaair))} \\&\leq C\|\zeta \|_{\rL^\infty(0,T;\rL^\infty(\T^2))}\| \mathrm{Z} \|_{\rL^\infty (0,T;\rL^{\infty}(\T^2))} \| V^\atm \|^2_{\rL^{2p}(0,T;\rH^{1+\nicefrac{1}{q},q}(\Omegaair))} \\& \leq C R^2.
        \end{aligned}
    \end{equation*}
where we used the embeddings
    \begin{equation}\label{eq:maxregembedding}
        \E^\atm_{1,T} \hookrightarrow \rL^{2p}(0,T; \rH^{2 -\nicefrac{1}{p},q}(\Omegaair)) \hookrightarrow \rL^{2p}(0,T;\rH^{1+\nicefrac{1}{q},q}(\Omegaair)),
      \end{equation}
which hold true whenever $\nicefrac{1}{p}+ \nicefrac{1}{q} \le 1$. Next, we estimate the nonlinear term including the vertical velocity anisotropically to obtain 
    \begin{equation*}
      \| \frac{\zeta}{\xi_0}W^\atm   \dz V^\atm \|_{\rLp(0,T;\rLq(\Omegaair))}  \leq C \|\zeta W^\atm \|_{\rL^\infty_t \rLq_\rH \rL^\infty_z} \| V^\atm \|_{\rLp_t \rL^\infty_\rH \rH^{1,q}_z}
      \leq C\left  \| \zeta W^\atm \right \|_{\rL^\infty_t \rLq_\rH \rL^\infty_z} \| V^\atm \|_{\E^\atm_{1,T}}.
    \end{equation*}
Combining the embeddings $\rH^{1,q}_z \hookrightarrow \rL^\infty_z$ , $\E^\atm_{1, T} \hookrightarrow \rL^\infty(0, T; \rL^\infty_\rH(G; \rLq_z(0, 1)))$ for $\nicefrac{1}{p}+\nicefrac{1}{q} < 1$  with 
the representation \eqref{eq:w}, we obtain 
    
    \begin{equation*}
        \begin{aligned}
             \left \| \zeta W^\atm \right \|_{\rL^\infty_t \rLq_\rH \rL^\infty_z} &\leq C \|  \widetilde{V}^\atm\mathrm{Z}^\top \nabla_\rH\zeta + \zeta \mathrm{Z}^\top : \nabla_\rH \widetilde{V}^\atm\|_{\rL^\infty (0,T;\rLq(\Omegaair))} \\
             &\leq C \| \widetilde{V}^\atm \|_{\rL^\infty (0,T;\rL^\infty_\rH \rLq_z)} \| \mathrm{Z} \|_{\rL^\infty (0,T;\rL^\infty(\T^2))} \| \nablaH \zeta \|_{\rL^\infty (0,T;\rLq(\T^2))} \\
             &\quad + C\| \zeta \|_{\rL^\infty (0,T;\rL^\infty(\T^2))}  \| \mathrm{Z} \|_{\rL^\infty (0,T;\rL^\infty(\T^2))} \| \nablaH \widetilde{V}^\atm \|_{\rL^\infty (0,T;\rLq(\Omegaair))} \\
             & \leq  C\| V \|_{\E^\atm_{1,T}}
        \end{aligned}
    \end{equation*}
    and therefore
    \begin{equation*}
          \| \frac{\zeta}{\xi_0}W^\atm   \dz V^\atm \|_{\rLp(0,T;\rLq(\Omegaair))}  \leq CR^2.
    \end{equation*}
Finally, we estimate the pressure terms and obtain first 
    \begin{equation*}
\begin{aligned}
         \big \|\frac{ \mathrm{Z}^\top \nablaH \zeta }{\xi_0 } \big\|_{\rLp(0,T;\rLq(\Omegaair))} \leq C \| \mathrm{Z} \|_{\rL^\infty (0,T;\rL^{\infty}(\T^2))}
      \| \zeta\|_{\rLp(0,T;\rH^{1,q}(\T^2))}\leq C T^{\nicefrac{1}{p'}}.
\end{aligned}
    \end{equation*}

The terms $F_1$ and $F_2$ are structurally the same and therefore the precise estimate of $F_2$ is omitted here.
    Concerning the nonlinear term $F_3$ we obtain
    \begin{equation}\label{eq: nonhom div}
        \| F_3 (V^\ocn) \|_{\mathbb{G}_{1,T}} \leq  \| \divH \bigr ( (\mathrm{J}\mathrm{Z}^\top-\mathbb{I}_2 ) \bar{V}^\ocn \bigl )\|_{\rH^{1,p}(0,T;\dot \rH^{-1,q}(\T^2))} + \| \nablaH \bar{V}^\ocn : \left [\mathrm{Z}^\top  - \mathbb{I}_2 \right ] \|_{\rLp(0,T;\rH^{1,q}(\T^2))}.
    \end{equation}
    Here, the second term can be estimated as follows
    \begin{equation*}
        \| \nablaH \bar{V}^\ocn : \left [\mathrm{Z}^\top  - \mathbb{I}_2 \right ] \|_{\rLp(0,T;\rH^{1,q}(\T^2))} \leq \| \nablaH \bar{V}^\ocn \|_{\rLp(0,T;\rH^{1,q}(\T^2))} \| \mathrm{Z}^\top - \mathbb{I}_2 \|_{\rL^\infty(0,T;\rH^{1,q}(\T^2))} \leq 
        CT^{\nicefrac{1}{p'}}.
    \end{equation*}
    Concerning the first term on the right-hand side of \eqref{eq: nonhom div}, it is sufficient to show that $\ (\mathrm{Z}-\mathbb{I}_2 ) \bar{V}^\ocn\in \rH^{1,p}(0,T;\rLq(\T^2))$, as the horizontal divergence maps $\rLq(\T^2,\R^2)$ to $\dot \rH^{-1,q}(\T^2;\R)$. We start by estimating the time derivative, i.\ e.,
    \begin{equation*}
        \begin{aligned}
            \| \dt (  (\mathrm{Z}-\mathbb{I}_2 ) \bar{V}^\ocn) \|_{\rLp(0,T;\rLq(\T^2))} &\leq \| \dt \mathrm{Z} \cdot  \nablaH \bar{V}^\ocn \|_{\rLp(0,T;\rLq(\T^2))} + \|  (\mathrm{Z}-\mathbb{I}_2 ) \dt \bar{V}^\ocn) \|_{\rLp(0,T;\rLq(\T^2))} \\
            & \leq C \| \dt \mathrm{Z} \|_{\rLp(0,T;\rH^{1,q}(\T^2))} \| V^\ocn \|_{\E^\ocn_{1,T}} + \|  (\mathrm{Z}-\mathbb{I}_2) \|_{\rL^\infty(0,T;\rL^\infty(\T^2))} \| V^\ocn \|_{\E^\ocn_{1,T}} \\ &\leq  CR^2+ C T^\delta ,
        \end{aligned}
    \end{equation*}
    where we used the identity
    \begin{equation*}
        \dt (\mathrm{Z} - \mathbb{I}_2) = \dt \mathrm{Z} = - (\nablaH \rX)^{-1} (\dt \nablaH \rX ) (\nablaH \rX)^{-1} = -\mathrm{Z} (\dt \nablaH \rX ) \mathrm{Z}
     = -\mathrm{Z} ( \nablaH \bar{V}^\air ) \mathrm{Z} 
    \end{equation*}
    see also \citep[Section~6]{HIRZ:25}. Since $F_4 \in \rLp(0,T;\rH^{1,q}(\T^2)) \hookrightarrow \rL^{p}(0,T;\rLq(\T^2))$ we conclude desired the estimate. Next, in order to estimate the nonlinear boundary term we embed $\rH^{1,q}(\T^2) \hookrightarrow \rB_{qq}^{1 - \nicefrac{1}{q}}(\T^2))$ and make use of the Banach algebra structure of the space $\rH^{1,q}(\T^2)$ to obtain
\begin{equation*}
    \begin{aligned}
        \|\zeta B(V^\atm,V^\ocn) \|_{\rLp(0,T;\rB_{qq}^{1 - \nicefrac{1}{q}}(\T^2))} &\leq C \|\zeta B(V^\atm,V^\ocn) \|_{\rLp(0,T;\rH^{1,q}(\T^2))} \\ 
        &\leq \|\zeta \|_{\rL^\infty(0,T;\rH^{1,q}(\T^2))} \| V^\air-V^\ocn \|^2_{\rL^{2p}(0,T;\rH^{1,q}(\T^2))}.
    \end{aligned}
\end{equation*}
In view of trace theorem and the embeddings \eqref{eq:maxregembedding} it follows that
\begin{equation*}
   \|\zeta B(V^\atm,V^\ocn) \|_{\rLp(0,T;\rB_{qq}^{1 - \nicefrac{1}{q}}(\T^2))}  \leq C  \|\zeta \|_{\rL^\infty(0,T;\rH^{1,q}(\T^2))} \bigl ( \| V^\air \|^2_{\E^\air_{1,T}} +  \| V^\ocn \|^2_{\E^\ocn_{1,T}} \bigr ) \leq C R^2.
\end{equation*}
Next, observe that
\begin{equation*}
    \rH^{\nicefrac{1}{2},p}(0,T;\rLq(\T^2)) \hookrightarrow \mathrm{F}^{\nicefrac{1}{2}}_{p2}(0,T;\rLq(\T^2)) \hookrightarrow \mathrm{F}^{\nicefrac{1}{2}- \nicefrac{1}{2p}}_{pq}(0,T;\rLq(\T^2)),
\end{equation*}
for $q \geq 2$. In view of $\rH^{\nicefrac{1}{2},p}(0,T)$ being a Banach algebra for $p>2$ and trace theorem we obtain

\begin{equation*}
    \begin{aligned}
        \| \zeta B(V^\atm, V^\ocn) \|_{\mathrm{F}_{pq}^{\nicefrac{1}{2}-\nicefrac{1}{2q}}(0,T;\rLq(\T^2))} &\leq C \|\zeta B(V^\atm,V^\ocn) \|_{\rH^{\nicefrac{1}{2},p}(0,T;\rLq(\T^2))} \\ 
        &\leq C  \| \zeta\|_{\rH^{\nicefrac{1}{2},p}(0,T; \rL^\infty( G))} \| V^\air-V^\ocn\|^2_{\rH^{\nicefrac{1}{2},p}(0,T;\rL^{2q}(\T^2))}\\
        &\leq  C \bigl ( \| V^\air \|^2_{\E^\air_{1,T}} +  \| V^\ocn \|^2_{\E^\ocn_{1,T}}  \bigr ) \\ &\leq CR^2.
    \end{aligned}
\end{equation*}
where we used the embeddings
\begin{equation*}
    \E^i_{1,T} \hookrightarrow \rH^{\nicefrac{1}{2},p}(0,T;\rH^{\nicefrac{2}{q},q}(\Omega^i)) \hookrightarrow \rH^{\nicefrac{1}{2},p}(0,T;\rW^{\nicefrac{1}{2q},2q}(\Omega^i)), 
\end{equation*}
which holds true whenever $q\geq 2$ and $i \in \{\air , \ocn \}$. 

\end{proof}
    
We are now in position to prove the local well-posedness result:

\begin{proof}[Proof of \autoref{thm:localwpCAO}]
We invoke a reference solution $( V^\atm_\ast, V^\ocn_\ast, \Pi_\ast) \in \E_{1,T}$ solving the system \eqref{eq:CAOatm} with homogeneous right-hand sides, suitable boundary conditions and initial value $(V^\atm_0,V^\ocn_0)$. The existence of reference solutions is guaranteed in view of \autoref{lem:Lame equation} and it holds true that $\| ( V^\atm_\ast, V^\ocn_\ast, \Pi_\ast) \|_{\E_{1,T}} \to 0$ as $T \to 0$. 
    With the prescript $_0$ indicating homogeneous initial data, we set
\begin{equation*}
    \begin{aligned}
          &\mathbb{B}_r(0,T)\coloneqq \{( V^\atm,  V^\ocn,\Pi) \in \E_{1,T} \text{ such that } \\ &\qquad (V^\atm- V^\atm_\ast, V^\ocn-V^\ocn_\ast) \in \, _0{\E^\atm_{1}}(0,T) \times \, _0{\E^\ocn_{1}}(0,T) \text{ and }  \| ( V^\atm-V^\atm_\ast, V^\ocn -V^\ocn_\ast,\Pi -\Pi_\ast ) \|_{\E_{1,T}}  \le r \}.
    \end{aligned}
\end{equation*}
In particular, we have
\begin{equation*}
    \| ( V^\atm,  V^\ocn,\Pi) \|_{\E_{1,T}} \leq r + \| ( V^\atm_\ast, V^\ocn_\ast, \Pi_\ast) \|_{\E_{1,T}}. 
\end{equation*}
For $( V^\atm_1, V^\ocn_1,\Pi_1 ) \in  \mathbb{B}_r(0,T)$, define the solution map $\Phi$ by 
\begin{equation*}
\Phi \colon \mathbb{B}_r(0,T) \to \E_{1,T}, \ ( V^\atm_1, V^\ocn_1,\Pi_1) \mapsto
    (V^\atm,V^\ocn, \Pi ),
\end{equation*}
where $(V^\atm,V^\ocn, \Pi )$ denotes the 
the solution to system \eqref{eq:CAOatm} with right-hand sides $F_1(\zeta_1, V_1^\atm)$, $F_2(V^\ocn_1, \Pi_1)$, $F_3(V_1^\ocn)$ as well as boundary conditions $\delta B(V^\atm_1,V^\ocn_1)$ and $\zeta_1 B(V^\atm_1,V^\ocn_1)$. Here $\zeta_1$ is defined by \eqref{eq:zetaJ}, that is 
\begin{equation*}
    \zeta_1=\frac{\xi_0}{\mathrm J_1},
\end{equation*}
with $\mathrm{J}_1$ denoting the flow component associated with $V^\air_1$.
By \autoref{lem:Lame equation} and \autoref{lem:nonlinearestimates}, $\Phi$ is well defined for $( V^\atm_1,  V^\ocn_1,\Pi_1 ) \in  \mathbb{B}_r(0,T)$ and we show next that $\Phi$ is a self-map and a contraction.
Notice that 
\begin{equation*}
     (V^\atm- V^\atm_\ast, V^\ocn-V^\ocn_\ast) \in \,_0{\E^\atm_{1}}(0,T) \times \, _0{\E^\ocn_{1}}(0,T)
\end{equation*}
and its left to show that
\begin{equation*}
    \| ( V^\atm-V^\atm_\ast, V^\ocn -V^\ocn_\ast,\Pi -\Pi_\ast ) \|_{\E_{1,T}}  \le r.
\end{equation*}
Using \autoref{lem:Lame equation} as well as the estimates \autoref{lem:nonlinearestimates} yields
\begin{equation*}
    \begin{aligned}
      \quad \| (V^\atm-V^\atm_\ast, V^\ocn -V^\ocn_\ast,\Pi -\Pi_\ast ) \|_{\E_{1,T}}  & \le C \bigl ( \| ( F_1(\zeta_1, V^\atm_1), F_2(V^\ocn_1, \Pi_1) , F_3(V^\ocn)) \|_{\E_{0,T}} \\& \quad + \| \delta B(V^\atm_1, V^\ocn_1) \|_{\F_{\N,T}} +   \| \zeta_1 B(V^\atm_1, V^\ocn_1) \|_{\F_{\N,T}} \bigr ) \\ & \leq C \big ( T^\delta + ( r + \| ( V^\atm_\ast, V^\ocn_\ast, \Pi_\ast) \|_{\E_{1,T}}  )^2\big) 
    \end{aligned}
\end{equation*}
where the constant $C > 0$ is independent of $T$ as $ (V^\atm-V^\atm_\ast, V^\ocn -V^\ocn_\ast,\Pi -\Pi_\ast )$ has vanishing initial data. 
Hence, choosing $T>0$ and $r>0$ sufficiently small we obtain 
 \begin{equation*}
     C \big ( T^\delta + ( r + \| ( V^\atm_\ast, V^\ocn_\ast, \Pi_\ast) \|_{\E_{1,T}}  )^2\big) \leq r
 \end{equation*}
 and therefore $\Psi$ is a maps $\mathbb{B}_r(0,T)$ to itself.
The contraction property can be shown similarly. 

Finally, choose $T$ sufficiently small such that $\rX(t,\cdot)$ is a $\rC^1$-diffeomorphism from $\T^2$ into $\T^2$ and set $\mathrm{Y}(t,\cdot) = \rX^{-1}(t,\cdot)$. Then defining 
\begin{equation*}
    \begin{aligned}
            \xi(t,x,y) &:= \zeta(t, \mathrm{Y}(t,x,y)), \ V^\atm(t,x,y,z) := V^\atm (t, \mathrm{Y}(t,x,y),z), \\   v^\ocn(t,x,y,z) &:= V^\ocn (t, \mathrm{Y}(t,x,y),z), \   \pi_s(t,x,y) := \Pi(t,\mathrm{Y}(t,x,y)).
    \end{aligned}
\end{equation*}
we verify that $(\xi, \vair, \vocn, \pi_s) \in \E_1(0,T)$ satisfies the system \eqref{eq:CAO1}. The uniqueness property follows from the uniqueness of the solution to the system \eqref{eq:CAOatm}.
\end{proof}

\section{Global in time existence}\label{sec:global}
 
\noindent
In this section, we show that the solution of system \eqref{eq:CAO1} remains arbitrarily close to the constant steady state $(\bar{\xi}, 0, 0, 0)$, provided the initial data is sufficiently close to it Here $\bar{\xi}$ denotes a positive constant. For this matter, we transform the atmospheric part of the system \eqref{eq:CAO1} by a slight modification of the Lagrange transform introduced in \autoref{subsec:lagrangecor}. More precisely we define the transformed atmospheric density by
the Jacobian of the flow via
\begin{equation}\label{eq:zetaJglobal}
    (\zeta+\bar\xi)\,\mathrm J = \xi_0,
    \ \text{ that is }\
    \zeta(t,y_\H)=\frac{\xi_0(y_\H)}{\mathrm J(t,y_\H)}-\bar\xi.
\end{equation}
The remaining quantities $(\vair,\vocn,\pi_s)$ are transformed as in \autoref{subsec:lagrangecor}.
Hence the transformed system takes the form
\begin{equation}
	\left\{
	\begin{aligned}
		\partial_t  V^\air
        - \frac{\mu \DeltaH V^\air }{(1-\delta z)\bar\xi}
        - \frac{\mu}{\bar\xi} \partial_z \bigl ( \frac{1-\delta z}{\delta^2} \partial_z V^\air \bigr )
        - \frac{\mu' \nablaH \divH V^\air}{(1-\delta z)\bar\xi}
        + \frac{\nablaH \zeta}{\bar\xi}
        &= F_1(\zeta,V^\air), &&\text{ in } \Omega^\air \times (0,T), \\
        \partial_t  V^\ocn - \Delta V^\ocn + \nablaH \Pi
        &= F_2(\Pi,V^\ocn), &&\text{ in } \Omega^\ocn \times (0,T), \\
		\divH  \bar{V}^{\ocn}
        &= F_3(V^\ocn), &&\text{ in } \T^2 \times (0,T),\\
        (V^\air(0),V^\ocn(0)) &= (v^\air_0,v^\ocn_0).
	\end{aligned}
	\right.
	\label{eq:CPEatmostransformed}
\end{equation}
where $\zeta$ is given by \eqref{eq:zetaJglobal}. The transformed boundary conditions read
\begin{equation}
\begin{aligned}\label{eq:nonlin bc}
     \partial_z V^\air|_{\Gamma^\air_u} = 0 \ &\text{ and } \ (\partial_z  V^\air)|_{\Gamma_i} = \delta (V^\ocn-V^\air)|V^\air - V^\ocn| , \\
   \partial_z V^\ocn|_{\Gamma^\ocn_b} = 0 \ &\text{ and } \ (\partial_z V^\ocn)|_{\Gamma_i} = (\zeta+ \bar{\xi}) (V^\air -V^\ocn)|V^\air - V^\ocn|. \\
\end{aligned}
\end{equation}

Here, the right-hand sides are given by
\begin{equation}
\label{eq:nonlinearities2}
\begin{aligned}
(F_1)_i
&=
\frac{\mu}{(1-\delta z)\bar\xi}
\Bigl(
\sum_{l,j,k}
\frac{\partial^2 V^\air_i}{\partial y_l \partial y_k}
\bigl(\mathrm Z_{k,j}-\delta_{k,j}\bigr)\mathrm Z_{l,j}
+
\sum_{l,k}
\frac{\partial^2 V^\air_i}{\partial y_l \partial y_k}
\bigl(\mathrm Z_{l,k}-\delta_{l,k}\bigr)
+
\sum_{l,j,k}
\mathrm Z_{l,j}
\frac{\partial V^\air_i}{\partial y_k}
\frac{\partial \mathrm Z_{k,j}}{\partial y_l}
\Bigr) \\
&\quad
+
\frac{\mu'}{(1-\delta z)\bar\xi}
\Bigl(
\sum_{l,j,k}
\frac{\partial^2 V^\air_j}{\partial y_l \partial y_k}
\bigl(\mathrm Z_{k,j}-\delta_{k,j}\bigr)\mathrm Z_{l,i}
+
\sum_{l,j}
\frac{\partial^2 V^\air_j}{\partial y_l \partial y_j}
\bigl(\mathrm Z_{l,i}-\delta_{l,i}\bigr)
+
\sum_{l,j,k}
\mathrm Z_{l,i}
\frac{\partial V^\air_j}{\partial y_k}
\frac{\partial \mathrm Z_{k,j}}{\partial y_l}
\Bigr) \\
&\quad
-
\frac{\zeta}{\bar\xi}(\partial_t V^\air)_i
-\frac{\zeta+\bar\xi}{\bar\xi}
\Bigl(
\sum_{l,j}
(V^\air-\bar V^\air)_l
\frac{\partial V^\air_i}{\partial y_j}\,
\mathrm Z_{l,j}
+
\frac{1-\delta z}{\delta}\,
W^\air \frac{\partial V^\air_i}{\partial z}
\Bigr)
-\frac{1}{\bar\xi}
\sum_j \frac{\partial \zeta}{\partial y_j}(\mathrm Z_{j,i}-\delta_{j,i}),
\\[0.3em]
(F_2)_i
&=
\sum_{l,j,k}
\frac{\partial^2 V^\ocn_i}{\partial y_l \partial y_k}
\bigl(\mathrm Z_{k,j}-\delta_{k,j}\bigr)\mathrm Z_{l,j}
+
\sum_{l,k}
\frac{\partial^2 V^\ocn_i}{\partial y_l \partial y_k}
\bigl(\mathrm Z_{l,k}-\delta_{l,k}\bigr)
+
\sum_{l,j,k}
\mathrm Z_{l,j}
\frac{\partial V^\ocn_i}{\partial y_k}
\frac{\partial \mathrm Z_{k,j}}{\partial y_l} \\
&\quad
-
\sum_{l,j}
(V^\ocn-\bar V^\air)_l
\frac{\partial V^\ocn_i}{\partial y_j}\,
\mathrm Z_{l,j}
-
W^\ocn \frac{\partial V^\ocn_i}{\partial z}
-
\sum_j
\frac{\partial \Pi}{\partial y_j}
\bigl(\mathrm Z_{j,i}-\delta_{j,i}\bigr),
\\[0.3em]
F_3
&=
-\nablaH \bar V^\ocn : \bigl[\mathrm Z^\top-\mathbb I_2\bigr].
\end{aligned}
\end{equation}

\subsection{Linearized CAO Problem and estimates of the nonlinear terms}
\label{subsec:Lincom}
\noindent
In this subsection we address the well-posedness of the linear system associated to \eqref{eq:CPEatmostransformed} and estimate the nonlinear terms given in \eqref{eq:nonlin bc} as well as \eqref{eq:nonlinearities2}. Concerning the oceanic part of \eqref{eq:CPEatmostransformed}, the linear system is the same as \eqref{eq:hstokes1} and well-posedness is guaranteed by \autoref{lem:Lame equation}(c). However, regarding the atmosphere, the system under considerations is now given by
\begin{equation}
	\left\{
	\begin{aligned}
		\dt  V^\air  - \frac{\mu \DeltaH V^\air }{(1-\delta z)}- \mu \dz \bigl ( \frac{1-\delta z}{\delta^2} \dz V^\air \bigr ) -  \frac{\mu' \nablaH \divH V^\air}{(1-\delta z)}  + \nablaH \zeta  &= f_1, &&\text{ in } \Omegaair \times (0,T),  \\
\dz        V^\air &=0, &&\text{ on }\T^2 \times \{1 \} \times (0,T), \\
        \dz V^\air &= b_1 , &&\text{ on } \T^2 \times \{0 \}  \times (0,T), \\
        V^\air(0)&= \vair_0,  
	\end{aligned}
	\right. 
	\label{eq:linearatmos}
\end{equation}
for a given time $0 < T < \infty$. Here, we have set $\bar{\xi}=1$ for simplicity of the notation.

\begin{prop}
    \label{prop: lin atmos}
     Let $T>0$, $p,q \in (2,\infty)$ and assume that the data $(f_1,b_1,\vair_0)$ satisfies \autoref{assu:data}. Then the system \eqref{eq:linearatmos} admits a unique, strong solution $V^\air \in  \E^\air_{1,T}.$
Furthermore, there is a constant $C>0$ such that
     \begin{equation}\label{eq: est atm}
         \begin{aligned}
             \| V^\air \|_{\E^\air_{1,T}} \leq C \bigl ( \| f_1 \|_{\rLp(0,T;\rLq(\Omegaair))} + \| b_1 \|_{\F_{\mathrm{N},T}}  + \| \vair_0 \|_{\rB^{2(1-\nicefrac{1}{p})}_{qp}(\Omegaair)} \bigr ).
         \end{aligned}
     \end{equation}
\end{prop}
\begin{proof}
    Solving system \eqref{eq:linearatmos} is equivalent to solving an inhomogeneous second order system 
    \begin{equation}
	\left\{
	\begin{aligned}
		\dt  V^\air  - \frac{\mu \DeltaH V^\air }{(1-\delta z)}- \mu \dz \bigl ( \frac{1-\delta z}{\delta^2} \dz V^\air \bigr ) -  \frac{\mu' \nablaH \divH V^\air}{(1-\delta z)}  &= f_1-\nablaH \zeta, &&\text{ in } \Omegaair \times (0,T),  \\
\dz        V^\air &=0, &&\text{ on } \T^2 \times \{1 \} \times (0,T), \\
        \dz V^\air &= b_1 , &&\text{ on } \T^2 \times \{0 \} \times (0,T), \\
        V^\air(0) &= \vair_0, 
	\end{aligned}
	\right. 
	\label{eq: linear Lame with pressure}
\end{equation}
where $\zeta$ is given by 
\begin{equation*}
    \zeta =\frac{\xi_0}{\mathrm{J}} - 1.
\end{equation*}
Consider the map 
\begin{equation*}
    R: \E^\air_{1,T} \to \E^\air_{1,T} \ \text{ with } \ V^\air \mapsto S_V (f_1- \nablaH \zeta, b_1, \vair_0),
\end{equation*}
where $S_V$ denotes the corresponding solution map of the equation \eqref{eq:Lame1}. By \autoref{lem:Lame equation} the map $S_V$ is continuous. The problem \eqref{eq: linear Lame with pressure} then admits a unique solution $ V^\air$ if and only if $R$ has a unique fixed point. 
First, we verify that $R : \E^\air_{1,T} \to \E^\air_{1,T}$ is well defined. This follows immediately from the fact that the data satisfies \autoref{assu:data} and using \autoref{lem:Lame equation}. In order to verify the contraction property of $R$, let $V_1^\air, V_2^\air \in \E_{1,T}^\air$ and 
\[
\mathrm J_i := \det \nablaH \mathrm X_i,
\qquad
\zeta_i := \frac{\xi_0}{\mathrm J_i}-1,
\qquad i\in\{1,2\},
\]
where $\mathrm X_i$ denotes the hydrostatic Lagrangian flow associated with $\bar V_i^\air$. Then
\[
R(V^\air_i)=S_V(-\nablaH \zeta_i,0,0),
\]
and therefore, by linearity of $S_V$,
\[
R(V^\air_1)-R(V^\air_2)
=
S_V(-\nablaH(\zeta_1-\zeta_2),0,0).
\]
It follows that
\begin{equation*}
    \| R(V^\air_1)- R(V^\air_2) \|
    \leq C \| \nablaH( \zeta_1 - \zeta_2) \|_{\rLp(0,T;\rLq(\T^2))}.
\end{equation*}
Moreover,
\[
\zeta_1-\zeta_2
=
\xi_0\Bigl(\frac{1}{\mathrm J_1}-\frac{1}{\mathrm J_2}\Bigr)
=
\xi_0\,\frac{\mathrm J_2-\mathrm J_1}{\mathrm J_1\mathrm J_2}.
\]
Since $\mathrm J_1$ and $\mathrm J_2$ are bounded away from zero and $\rH^{1,q}(\T^2)$ is a Banach algebra for $q>2$, we obtain
\begin{equation*}
    \| \zeta_1-\zeta_2 \|_{\rL^\infty(0,T;\rH^{1,q}(\T^2))}
    \leq
    C
    \| \mathrm J_1-\mathrm J_2 \|_{\rL^\infty(0,T;\rH^{1,q}(\T^2))}.
\end{equation*}
Using H\"older's inequality, we infer
\begin{equation*}
    \| \nablaH(\zeta_1-\zeta_2) \|_{\rLp(0,T;\rLq(\T^2))}
    \leq
    C T^{\nicefrac{1}{p}}
    \| \zeta_1-\zeta_2 \|_{\rL^\infty(0,T;\rH^{1,q}(\T^2))}
    \leq
    C T^{\nicefrac{1}{p}}
    \| \mathrm J_1-\mathrm J_2 \|_{\rL^\infty(0,T;\rH^{1,q}(\T^2))}.
\end{equation*}
Finally, since $\mathrm J_1(0)=\mathrm J_2(0)=1$, the difference $\mathrm J_1-\mathrm J_2$ has vanishing initial trace, and the Lipschitz dependence of the Jacobian on the flow yields
\begin{equation*}
    \| \mathrm J_1-\mathrm J_2 \|_{\rL^\infty(0,T;\rH^{1,q}(\T^2))}
    \leq
    C
    \| V^\air_1 -V^\air_2 \|_{\E^\air_{1,T}},
\end{equation*}
with $C>0$ independent of $T$. Consequently,
\begin{equation*}
     \| R(V^\air_1)- R(V^\air_2) \|
     \leq
     C T^{\nicefrac{1}{p}}
     \| V^\air_1 -V^\air_2 \|_{\E^\air_{1,T}}
     \leq
     \kappa  \| V^\air_1 -V^\air_2 \|_{\E^\air_{1,T}},
\end{equation*}
with $\kappa \in (0,1)$ for $T>0$ sufficiently small.The contraction mapping principle then yields the existence of a unique $V^\air_* \in \E^\air_{1,T}$ such that $R V^\air_* = V^\air_*$, i.\ e., a unique solution to \eqref{eq: linear Lame with pressure}.
To extend the above result to arbitrary $T>0$ we employ the following iteration scheme. Let $\frac{9}{10} \cdot T' < T_1 < T'$, where $T'$ is the largest $T$ such that $R$ is a contraction mapping, and $V^{\air,1} \in \E_{1,T_1}^\air \hookrightarrow \mathrm{BUC}([0,T_1];\rB_{qp,\per}^{2(1-\nicefrac{1}{p})}(\Omega^\air;\R^2))$ the solution of \eqref{eq: linear Lame with pressure}. 
        We consider the problem \eqref{eq: linear Lame with pressure} on the time interval $[T_1,2T_1]$ with initial data $V^\air(T_1) = V^{\air,1}(T_1)$. 
		Repeating the argument above, we obtain the existence of a unique solution $V^{\air,2}$ of \eqref{eq: linear Lame with pressure} on $[T_1,2T_1]$. 
        Since $[0,T'] \cap [T_1,2T_1] \neq \emptyset$, we deduce from the uniqueness that $V^{\air,1}|_{[0,T'] \cap [T_1,2T_1]} = V^{2,\air}|_{[0,T'] \cap [T_1,2T_1]}$.
        Note that on the interval \(I_k=[T_k,T_{k+1}]\) we keep the Lagrangian flow
constructed on the preceding intervals. Thus, for all candidates in the
fixed point argument on \(I_k\), the initial flow is fixed by
\[
\rX_i(T_k)=\rX_\ast(T_k),\quad \mathrm J_i(T_k)=\mathrm J_\ast(T_k),\quad i=1,2,
\]
where $\rX_\ast$ denotes the previous flow on $I_{k-1}$.
Consequently,
\[
(\mathrm J_1- \mathrm J_2)(T_k)=0.
\]
Hence the same Lipschitz estimate for \( \mathrm J_1- \mathrm J_2\) as on the first interval
applies, with \(0\) replaced by \(T_k\). 
        As a result, we get a unique solution $V \in \E_{1,T}^\air(0,2T_1)$ of \eqref{eq: linear Lame with pressure}. Iterating this argument leads to a unique solution $V \in \E_{1,T}^{\air}$ and the a-priori estimate is a consequence of the closed graph theorem.
\end{proof}

Next, we show that the hydrostatic Lagrangian transform remains well-defined for all large times $T > 0$, provided that the atmospheric velocity remains sufficiently close to the equilibrium state. For this purpose, consider $T>0$, $p,q \in (2,\infty)$ and $\rX$ as in \autoref{subsec:lagrangecor}. Define the solution ball
\begin{equation}\label{eq: solball}
    \mathbb{B}_\varepsilon := \{ ( V^\air, V^\ocn, \Pi) \in \E_{1,T} : \|(V^\air, V^\ocn, \Pi) \|_{\E_{1,T}} \leq \varepsilon \},
\end{equation}
with 
\begin{equation*}
    \| (V^\air, V^\ocn, \Pi) \|_{\E_{1,T}} :=  \| V^\air \|_{\E^\air_{1,T}} + \| V^\ocn \|_{\E^\ocn_{1,T}} + \| \Pi \|_{\rLp(0,T;\H^{1,q}(\T^2))}.
\end{equation*}

Then there is a constant $C>0$ such that for all $\varepsilon \in ( 0, \frac{1}{2CT^{1/p'}})$ we have 
\begin{equation*}
    \begin{aligned}
        \| \nablaH \rX - \mathbb{I}_2\|_{\rL^\infty((0,T) \times G)}  \leq C \int_0^T  \| \nablaH \bar{V}^\air(t,\cdot) \|_{\rH^{1,q}(\T^2)} \d t \leq CT^{\nicefrac{1}{p'}}\| V^\air \|_{\E^\air_{1,T}} \leq \frac{1}{2}.
    \end{aligned}
\end{equation*}
We conclude that choosing $\eps>0$ sufficiently small implies that the transform $\rX$ is a well-defined $\rC^1$-diffeomorphism. 
\begin{lem}\label{lem:ests of trafo global wp}
Let $T>0$, $p,q \in (2,\infty)$ and consider $(\zeta, V^\air, V^\ocn, \Pi) \in \mathbb{B}_{\eps}$, where $\mathbb{B}_{\eps}$ has been made precise in \eqref{eq: solball}.Let $\rX$, $\mathrm{Z}$ and $\mathrm{J}$ be the flow quantities associated with the averaged velocity $\bar V^\air$
Then the following statements hold true.

\begin{enumerate}[(a)]
    \item There is a constant $C>0$, independent of $\eps$, such that
    \begin{equation*}
        \sup_{t \in (0,T)} \|\nablaH \rX(t,\cdot)-\mathbb I_2\|_{\rL^\infty(\T^2)}
        + \sup_{t \in (0,T)} \|\nablaH \rX(t,\cdot)-\mathbb I_2\|_{\rH^{1,q}(\T^2)}
        \le C\eps.
    \end{equation*}
    If $\eps \in \bigl(0,\frac{1}{2C}\bigr)$, then in particular
    \[
        \|\nablaH \rX(t,\cdot)-\mathbb I_2\|_{\rL^\infty(\T^2)} \le \frac12
        \qquad\text{for all } t\in(0,T),
    \]
    so that $\nablaH \rX(t,\cdot)$ is invertible and $\mathrm Z=[\nablaH \rX]^{-1}$ is well defined.

    \item There is a constant $C>0$, independent of $\eps$, such that
    \begin{equation*}
        \|\mathrm Z\|_{\rL^\infty(0,T;\rH^{1,q}(\T^2))}
        + \|\mathrm Z^\top\|_{\rL^\infty(0,T;\rH^{1,q}(\T^2))}
        + \|\nablaH \rX\|_{\rL^\infty(0,T;\rH^{1,q}(\T^2))}
        \le C.
    \end{equation*}

    \item Similar estimates as in~(a) are also valid for $\mathrm Z$ as well as $\mathrm Z^\top$, i.\ e.,
    \begin{equation*}
        \sup_{t\in(0,T)} \|\mathrm Z(t,\cdot)-\mathbb I_2\|_{\rH^{1,q}(\T^2)}
        + \sup_{t\in(0,T)} \|\mathrm Z^\top(t,\cdot)-\mathbb I_2\|_{\rH^{1,q}(\T^2)}
        \le C\eps
    \end{equation*}
    for a constant $\eps>0$ independent of $T$.

    \item It holds that $\mathrm J(t,y_\H)\ge c>0$ for all $t\in(0,T)$ and $y_\H\in\T^2$. In particular, $\mathrm J^{-1}$ is well defined. Moreover, there is a constant $C>0$, independent of $\eps$, such that
  \begin{equation*}
        \| \mathrm J \|_{\rH^{1,p}(0,T;\rH^{1,q}(\T^2))}
        + \| \mathrm J \|_{\rL^\infty(0,T;\rH^{1,q}(\T^2))}
        \le C \ \text{ and } \  \sup_{t\in(0,T)} \| \mathrm J(t,\cdot)-1 \|_{\rH^{1,q}(\T^2)}
        \le C\eps.
    \end{equation*}
    \item For all $j,k,l \in \{1,2\}$, there is a constant $C>0$, independent of $T$, such that
    \begin{equation*}
        \Bigl\|\frac{\partial \mathrm Z_{l,j}}{\partial y_k}\Bigr\|_{\rL^\infty(0,T;\rL^q(\T^2))}
        \le C\eps.
    \end{equation*}
\end{enumerate}
\end{lem}
We are now in position to estimate the nonlinear terms as made precise in \eqref{eq:nonlin bc} and \eqref{eq:nonlinearities2}. 

\begin{prop}
    \label{prop:estnonlin}   
Let  $T>0$, $p,q \in (2,\infty)$ and $\eps>0$ sufficiently small so that the transformation above is well defined. Assume that $( V^\air, V^\ocn,\Pi)  \in \mathbb{B}_{\varepsilon}$
and $\| \xi_0 -1 \|_{\rH^{1,q}(\Omegaair)} \leq \eps/2$.
Then there is a constant $C>0$ depending on $p$, $q$ and $T$, it holds true that
    \begin{equation*}
        \| (F_1,F_2,F_3) \|_{\E_{0,T}}  + \|(B_1,B_2) \|_{\F^2_{\rN,T}}\leq C \varepsilon^2.
    \end{equation*}
\end{prop}
\begin{proof} 
First we adapt the estimates of the density. By \eqref{eq:zetaJglobal}, we have
\[
\zeta=\frac{\xi_0}{\mathrm J}-\bar\xi
=
\frac{\xi_0-\bar\xi}{\mathrm J}
+
\bar\xi\Bigl(\frac1{\mathrm J}-1\Bigr).
\]
Since $\mathrm J$ is bounded away from zero by Lemma~\ref{lem:ests of trafo global wp}(d), it follows from the Banach algebra property of $\rH^{1,q}(\T^2)$ for $q>2$, the estimate for $\mathrm J-1$ from Lemma~\ref{lem:ests of trafo global wp}(d), and standard composition results that
\begin{equation}\label{eq:zetaestimateglobal}
    \|\zeta\|_{\rL^{\infty}(0,T;\rH^{1,q}(\T^2))}
    \le C\eps,
    \qquad
    \|\zeta\|_{\rLp(0,T;\rH^{1,q}(\T^2))}
    \le C \eps,
\end{equation}
where $C>0$ is independent of $\eps$. From here on the proof uses the same techniques as the proof of \autoref{lem:nonlinearestimates}.
For convenience, this is demonstrated in detail for one term of $F_1$. The Banach algebra structure of $\rH^{1,q}(\T^2)$ yields
    \begin{equation*}
        \begin{aligned}
          &\big \|  \sum_{l,j,k} \frac{\partial^2 V^\atm_i}{\partial y_l \partial y_k} (\mathrm{Z}_{k,j} - \delta_{k,j}) \mathrm{Z}_{l,j} \big \|_{\rLp(0,T;\rLq(\Omegaair))}
          \\&\leq C \| V^\atm \|_{\rLp(0,T;\rH^{2,q}(\Omegaair))} \| \mathrm{Z} - \mathbb{I}_2 \|_{\rL^\infty (0,T;\rH^{1,q}(\T^2))} \| \mathrm{Z} \|_{\rL^\infty (0,T;\rH^{1,q}(\T^2))} \\
          &\leq C\eps^2.
        \end{aligned}
     \end{equation*}   
\end{proof}

We are now in position to prove \autoref{thm:globalWP}

\begin{proof}[Proof of \autoref{thm:globalWP}]
    Let $\mathbb{B}_\varepsilon$ be given as in \eqref{eq: solball}. For $( V^\air_1, V^\ocn_1,\Pi_1) \in  \mathbb{B}_\varepsilon$, we denote by 
    \begin{equation*}
        (V^\air, V^\ocn,\Pi) \coloneqq \Psi( V^\air_1, V^\ocn_1,\Pi_1) \in\E_{1,T}
    \end{equation*}
    the unique solution to the linearized equations \eqref{eq:CPEatmostransformed} with right-hand sides $F_1(\zeta_1^\air, V_1^\air)$, $F_2(\Pi_1, V_1^\ocn)$, $F_3(V^\ocn_1)$ as well as boundary conditions $B_1(V_1^\air,V_1^\ocn)$ and $B_2(\zeta_1, V_1^\air, V_1^\ocn)$ guaranteed by \autoref{prop:estnonlin}, where $\zeta_1$ is given by 
    \begin{equation*}
        \zeta_1 = \frac{\xi_0}{\mathrm{J}_1} - 1.
    \end{equation*}
    We will show that $\Psi$ is a self-map and a contraction. Choosing the initial data small enough, i.\ e.,
\begin{equation*}
    \| \xi_0-\bar{\xi} \|_{\rH^{1,q}(\T^2)}\leq \frac{\eps}{2}, \quad \| v^\air_0 \|_{\rB^{2(1-\nicefrac{1}{p})}_{qp}(\Omegaair)} \leq \frac{\varepsilon}{2C_1} \ \text{ and } \ \| \vocn_0 \|_{\rB^{2(1-\nicefrac{1}{p})}_{qp}(\Omegaocn)}  \leq \frac{\varepsilon}{2C_2},
\end{equation*}
where $C_1, C_2>0$ denote the constants given in the estimates of \autoref{prop: lin atmos} and \autoref{lem:Lame equation}(c) respectively. Using \autoref{prop: lin atmos}, \autoref{lem:Lame equation} and \autoref{prop:estnonlin} yields
 \begin{equation*}
        \begin{aligned}
            \|  ( V^\air, V^\ocn, \Pi) \|_{\E_{1,T} } 
            &\leq C \bigl ( \| v^\air_0 \|_{\rB^{2(1-\nicefrac{1}{p})}_{qp}(\Omegaair)} + \| \vocn_0 \|_{\rB^{2(1-\nicefrac{1}{p})}_{qp}(\Omegaocn)} \\& \quad 
            + \| (F_1,F_2,F_3) \|_{\E_{0,T}}+ \| F_4 \|_{\mathbb{G}_{1,T}}  + \|(B_1,B_2) \|_{\F^2_{\rN,T}}   \bigr ) \\
            &\leq \frac{\varepsilon}{2C_1}+\frac{\varepsilon}{2C_2}+C\varepsilon^2.
        \end{aligned}
    \end{equation*}
We conclude that for $\varepsilon>0$ small enough, $\Psi$ is a self-map. The contraction property follows in the same way. Finally, since $\varepsilon$ is chosen to be sufficiently small, the transform $\rX(t,\cdot )$ is a $\rC^1$-diffeomorphism from $\T^2$ into $\T^2$ and we set $\mathrm{Y}(t,\cdot) = \rX^{-1}(t,\cdot)$. Defining 
\begin{equation*}
    \begin{aligned}
            \xi(t,x,y) :=1+ \zeta(t, \mathrm{Y}(t,x,y)), \quad v^\air(t,x,y,z) &:= V^\air (t, \mathrm{Y}(t,x,y),z), \quad   v^\ocn(t,x,y,z) := V^\ocn (t, \mathrm{Y}(t,x,y),z), \\  \pi_s(t,x,y) &:= \Pi(t,\mathrm{Y}(t,x,y)).
    \end{aligned}
\end{equation*}
we verify that $(\xi, \vair,\vocn, \pi_s)\in \E_{1,T}$ is a global strong solution of \eqref{eq:CAO1}.
\end{proof}

\end{document}